\documentclass[12pt
]{article}

\newcommand{\noop}[1]{}

\usepackage{amssymb}
\usepackage{amsmath}
\usepackage{mathrsfs}
\usepackage{theorem}
\def\math{\mathscr}
\usepackage{newtxtext}

\newtheorem{theorem}{Theorem}
\newtheorem{lemma}{Lemma}
\newtheorem{corollary}{Corollary}

{\theorembodyfont{\rmfamily} }
{\theorembodyfont{\rmfamily} }
{\theorembodyfont{\rmfamily} \newtheorem{example}{Example}}

\newcommand{\C}{\mathbb{C}}

\newcommand{\n}{\nu}

\newcommand{\D}{\Delta}

\newcommand{\de}{\partial}
\newcommand{\ga}{\gamma}

\newcommand{\vf}{\varphi}

\newcommand{\si}{\sigma}

\newcommand{\GA}{\Gamma}
\newcommand{\al}{\alpha}

\newcommand{\la}{\lambda}

\newcommand{\ze}{\zeta}

\newcommand{\Om}{\Omega}
\newcommand{\om}{\omega}
\newcommand{\ro}{\varrho}
\newcommand{\dive}{\mathop{\rm div}\nolimits}

\newcommand{\I}{{\math I}}

\newcommand{\lan}{\langle}
\newcommand{\ran}{\rangle}

\renewcommand\leq{\leqslant}
\renewcommand\geq{\geqslant}

\newcommand{\R}{\mathbb R}
\newcommand{\CC}{\mathbb{C}}

\newcommand\Cspt{C_{0}}
\newcommand{\Hspt}{\mathaccent"017{H}}

\renewcommand\Re{\mathop{\mathbb R \rm{e}}\nolimits}
\renewcommand\Im{\mathop{\mathbb I \rm{m}}\nolimits}
\newcommand\elle{\mathop{\mathscr L}\nolimits}
\newcommand\A{\mathop{\mathscr A}\nolimits}
\newcommand\B{\mathop{\mathscr B}\nolimits}

\newcommand\Dm{\mathop{\mathscr D}\nolimits}

\def\Xint#1{\mathchoice 
  {\XXint\displaystyle\textstyle{#1}}%
  {\XXint\textstyle\scriptstyle{#1}}%
  {\XXint\scriptstyle\scriptscriptstyle{#1}}%
  {\XXint\scriptscriptstyle\scriptscriptstyle{#1}}%
  \!\int}
\def\XXint#1#2#3{{\setbox0=\hbox{$#1{#2#3}{\int}$} 
  \vcenter{\hbox{$#2#3$}}\kern-.5\wd0}}
\def\dashint{\Xint-}

\title{A survey of functional and $L^p$ dissipativity theory}

\author{A. Cialdea
\thanks{Department of Mathematics, Computer Sciences and Economics,
University of Ba\-si\-li\-ca\-ta, V.le dell'Ateneo Lucano, 10, 85100 Potenza, Italy.
 \textit{email:}
cialdea@email.it.}\and
V. Maz'ya 
\thanks{Department of Mathematics, Link\"oping University,
SE-581 83, Link\"oping, Sweden.
Peoples' Friendship University
 of Russia (RUDN University);
6 Miklukho-Maklaya St, Moscow, 117198, Russian Federation.
\textit{email}: vladimir.mazya@liu.se.}
}
\date{}    

\begin{document}

\maketitle

\bigskip\bigskip

 {\small \textbf{Abstract.} Various notions of dissipativity type for  partial differential operators 
 and their applications are surveyed. We deal with functional dissipativity and its particular case $L^p$-dissipativity. Most of the results are due to the authors.}

 \section{Introduction}

 The present paper contains a survey of recent results concerning dissipativity of partial differential operators.
 To be more precise, we mean the notion of functional dissipativity introduced in 
 \cite{CM2021} and its particular case, the  so called $L^p$-dissipativity.
 
Our joint studies in this area started in 2005, when we found
 necessary and sufficient conditions for  the
$L^p$-dissipativity of second order differential operators  with complex valued coefficients.

The $L^p$-dissipativity of a partial differential operator arises in a natural
 way in the study of partial differential equations with data in $L^p$. 
The theory of such problems has a long history.
 In fact  $L^p$-dissipativity appeared in 1937 in the pioneering work of Cimmino
   \cite{cimmino} on the Dirichlet problem with boundary data in $L^p$. Similar ideas were
   used in  \cite{mazya59} and \cite{mazyasobolev}.
In  \cite{Moblique}  the study of degenerate oblique derivative problem hinges
on the weighted $L^p$ positivity of the differential operator. 
Later we give more historical information.

In order to introduce the topic in a simple way, let us consider the classical Cauchy-Dirichlet problem for
the heat equation
\begin{equation}
    \begin{cases}
    \displaystyle\frac{\de u}{\de t}=\D u, & \text{for $t>0$},\\
    u(x,0)=\vf(x), & \text{$x\in\R^n$,}
    \end{cases}
    \label{eq:intr1}
\end{equation}
where $\vf$ is a given function in $C^{0}(\R^n)\cap 
L^{\infty}(\R^n)$.

It is well known that the unique solution of problem \eqref{eq:intr1}
in the class of smooth bounded solutions
is given by the formula
\begin{equation}
    u(x,t)=\frac{1}{\sqrt{(4\pi t)^{n}}}\int_{\R^n}\vf(y)\,
            e^{-\frac{|x-y|^{2}}{4t}}dy,\quad
	    x\in\R^n, t>0.
    \label{eq:sol}
\end{equation}

From \eqref{eq:sol} it follows immediately
\begin{equation}
    |u(x,t)|\leq \Vert \vf\Vert_{\infty},
    \quad t>0,
    \label{eq:diseq1}
\end{equation}
since
\begin{equation}
   \frac{1}{\sqrt{(4\pi t)^{n}}}\int_{\R^n}
                e^{-\frac{|x-y|^{2}}{4t}}dy=1
            \quad (t>0).
    \label{eq:=1}
\end{equation}

Inequality \eqref{eq:diseq1} leads to the classical maximum modulus principle
$$
    \Vert u(\cdot,t)\Vert_{\infty} \leq \Vert \vf \Vert_{\infty},
    \quad t>0,
$$
and this in turn implies that the norm $\Vert 
u(\cdot,t)\Vert_{\infty}$ is a decreasing function of $t$.
In fact,  fix $t_{0}>0$ and consider the problem
\begin{equation}
    \begin{cases}
     \displaystyle\frac{\de v}{\de t}=\D v, & \text{for $t>t_{0}$},\\
    v(x,t_{0})=u(x,t_{0}), & \text{$x\in\R^n$.}
    \end{cases}
    \label{eq:intr2}
\end{equation}

It is clear that the unique solution of \eqref{eq:intr2} is given by
$v(x,t)=u(x,t)$ ($t>t_{0}$) and we have
$$
\Vert v(\cdot,t)\Vert_{\infty}\leq \Vert u(\cdot,t_{0})\Vert_{\infty},
\quad t> t_{0},
$$
i.e.,
$$
\Vert u(\cdot,t)\Vert_{\infty}\leq \Vert u(\cdot,t_{0})\Vert_{\infty},
\quad t> t_{0}.
$$

The $L^{\infty}$ norm is not the only norm for which we
have this kind of dissipativity. Let us consider the $L^{p}$-norm
with $1<p<\infty$. By Cauchy-H\"older inequality, from \eqref{eq:sol} 
we get
$$
|u(x,t)|\leq \left(
\frac{1}{\sqrt{(4\pi t)^{n}}}\int_{\R^n}|\vf(y)|^{p}
            e^{-\frac{|x-y|^{2}}{4t}}dy\right)^{\!\!1/p}
	   \!\! \left(\frac{1}{\sqrt{(4\pi t)^{n}}}\int_{\R^n}
            e^{-\frac{|x-y|^{2}}{4t}}dy\right)^{\!\!1/p'}
$$
($1/p+1/p'=1$) and then, keeping in mind \eqref{eq:=1},
$$
|u(x,t)|^{p}\leq \frac{1}{\sqrt{(4\pi t)^{n}}}\int_{\R^n}|\vf(y)|^{p}
            e^{-\frac{|x-y|^{2}}{4t}}dy.
$$

Integrating over $\R^n$ and applying Tonelli's Theorem we find
\begin{gather*}
\int_{\R^n}|u(x,t)|^{p}dx \leq \frac{1}{\sqrt{(4\pi t)^{n}}}\int_{\R^n}dx\int_{\R^n}|\vf(y)|^{p}
            e^{-\frac{|x-y|^{2}}{4t}}dy=\\
\frac{1}{\sqrt{(4\pi t)^{n}}}\int_{\R^n}|\vf(y)|^{p}dy\int_{\R^n}
            e^{-\frac{|x-y|^{2}}{4t}}dx= 	
	    \int_{\R^n}|\vf(y)|^{p}dy  
\end{gather*}
and we have proved that
\begin{equation}
    \Vert u(\cdot,t)\Vert_{p}\leq \Vert \vf\Vert_{p}.
    \label{eq:dissip}
\end{equation}

As before, this inequality implies that the norm $\Vert 
u(\cdot,t)\Vert_{p}$ is a decreasing function of $t$.

Let us consider now the more general problem
\begin{equation}
    \begin{cases} 
    \displaystyle\frac{\de u}{\de t}=A u, & \text{for $t>0$},\\
    u(x,t)=0, & \text{ for $x\in\de\Om,\ t>0$},\\
    u(x,0)=\vf(x), & \text{$x\in\Om$,}
    \end{cases}
    \label{eq:intrE}
\end{equation}
where $\Om$ is a domain in $\R^{n}$ and $A$ is a linear
 elliptic partial differential operator of the second order
$$
    Au=\sum_{|\al|\leq 2} a_{\al}(x)\, D^{\al}u\, ,
$$
the coefficients $a_{\al}$ being complex valued. 
A natural question arises: under which conditions for the
operator $A$ the solution $u(x,t)$ of the problem \eqref{eq:intrE} 
satisfies the inequality \eqref{eq:dissip} ?

 As we 
know already, \eqref{eq:dissip} implies that $\Vert u(\cdot, 
t)\Vert_{p}$ is a decreasing function of $t$ and then
\begin{equation}
    \frac{d}{dt} \Vert u(\cdot, 
t)\Vert_{p}\leq 0.
    \label{eq:deriv<0}
\end{equation}
On the other hand, at least formally, we have for $1<p<\infty$,
\footnote{Note that $\de_{t}|u|=\de_{t}\sqrt{u\, \overline{u}}
=(u_{t}\overline{u}+u\overline{u}_{t})/(2\sqrt{u\,\overline{u}})=
\Re (u_{t}\overline{u}/|u|)$.}
\begin{equation}
    \frac{d}{dt} \Vert u(\cdot, 
    t)\Vert_{p}^{\ p}= \frac{d}{dt} \int_{\Om}|u(x,t)|^{p}dx =
    p \Re \int_{\Om}\lan \de_{t}u,u\ran |u|^{p-2}dx,
    \label{eq:dt}
\end{equation}
where $\lan \cdot, \cdot\ran$ denotes the usual scalar product in $\C$.
Since $u$ is a solution of the problem \eqref{eq:intrE}, 
keeping in mind  \eqref{eq:dt}, we have that \eqref{eq:deriv<0} holds
if, and only if,
$$
    \Re \int_{\Om}\lan Au,u\ran |u|^{p-2}dx \leq 0.
$$
This leads to the following definition: \textit{let $A$ a linear operator
from $D(A)\subset L^{p}(\Om)$ to $L^{p}(\Om)$; $A$ is said to be 
$L^{p}$-dissipative if }
\begin{equation}
    \Re \int_{\Om}\lan Au,u\ran |u|^{p-2}dx \leq 0, \quad
    \forall\ u\in D(A).
    \label{eq:Lpdiss}
\end{equation}
It follows from what we have said before that if $A$ is $L^{p}$-dissipative and if 
 problem 
\eqref{eq:intrE} has a solution, then \eqref{eq:deriv<0} holds.
Here we shall not dwell upon details of rigourous justification of the above argument.

We conclude Introduction by a well known result
(see e.g. \cite[p.215]{pazy}).
Consider the operator in divergence form
with real smooth coefficients 
\begin{equation}\label{eq:defApazy}
Au =  {\partial_{i}}\left(a_{ij}(x)\, \partial_j u
\right)
\end{equation}
($a_{ji}=a_{ij}\in 
C^{1}(\overline{\Om})$):
if  $a_{ij}(x)\xi_{i}\xi_{j}\geq 0$
for any $\xi\in\R^{n}$, $x\in \Om$, the operator \eqref{eq:defApazy} is 
 $L^{p}$-dissipative for any $p$.
  If $2\leq p<\infty$
this can be deduced easily by integration by parts. 
Indeed we
can write
\begin{gather*}
\int_{\Om} \langle Au , u\rangle |u|^{p-2}dx= -
\int_{\Om} a_{ij}\ \partial_j u\, 
 \partial_{i} \left(\overline{u}\, |u|^{p-2}\right)dx \\
= -
\int_{\Om} a_{ij}  \partial_j u  \left(|u|^{p-2} \partial_i \overline{u}+\overline{u}
\, {\partial_{i}}(|u|^{p-2})\right) d x \, .
\end{gather*}

Since
$$
\partial_{i} (|u|^{p-2})=
(p-2)|u|^{p-3} \partial_{i} |u| =
(p-2)|u|^{p-4}\Re \left(\overline{u}\, \partial_i u \right)
$$
we can write
\begin{gather*}
\partial_{j} u  \left(|u|^{p-2}\partial_i \overline{u}+\overline{u}\,
\partial_{i}(|u|^{p-2})\right)\\
=
|u|^{p-4}\left[\overline{u}\, u\, \partial_{j} u \,  \partial_i \overline{u}+(p-2) \overline{u}\, \partial_{j} u \Re\left( u\, \partial_i \overline{u}\right)\right].
\end{gather*}

Setting
$$
|u|^{(p-4) / 2} \overline{u}\, \partial_{j} u=\xi_{j}+i \eta_{j} \, ,
$$
we have
$$
a_{ij} \partial_{j} u  \left(|u|^{p-2}\partial_i \overline{u}+\overline{u}\,
\partial_{i}(|u|^{p-2})\right)
= a_{i j}\left(\xi_{j}+i \eta_{j}\right) 
{\left(\xi_{i}-i \eta_{i}+(p-2) \xi_{i}\right)} .
$$

This implies
$$
\Re\left(
a_{ij} \partial_{j} u  \left(|u|^{p-2}\partial_i \overline{u}+\overline{u}\,
\partial_{i}(|u|^{p-2})\right)\right)
= (p-1) a_{i j} \xi_{i} \xi_{j} 
+a_{i j} \eta_{i} \eta_{j}
$$
and then
$$
\Re \int_{\Om} a_{ij} \partial_{j} u\, 
 \partial_{i}\left(\overline{u}\, |u|^{p-2}\right)dx \geq 
0 ,
$$
i.e., $A$ is $L^{p}$-dissipative. Some extra arguments
are necessary for the case $1< p <2$.


During the last half a century various aspects of the
$L^{p}$-theory of semigroups generated by linear differential operators
were studied in
\cite{brezis, davies, amann, strichartz2,
davies1, kovalenko, robinson, davies2, davies3,
liskevich1, liskevichsemen, langer,  auscher, 
daners, karrmann, liskevich2, sobol, ouha, CM2005, 
metafune,CM2006} and others.
A general account of the subject can be found
in the book \cite{ouhabook}. Certain  of our earlier results have been described in our monograph \cite{CMbook},  where they are considered in the more general frame of semi-bounded operators.

Necessary and sufficient conditions for the
$L^{\infty}$-contractivity for general second order
strongly elliptic systems
with smooth coefficients
were given in \cite{ KM1994} (see also the monograph
\cite{kresinmazbook}).
Scalar second order elliptic operators with
complex coefficients were handled as a particular case.
The operators generating
$L^{\infty}$-contractive semigroups were later
characterized in \cite{auscher}
under the assumption that the coefficients are measurable and bounded.

The  maximum modulus principle for linear elliptic equations and systems 
with complex coefficients was considered 
by Kresin and Maz'ya. 
They have obtained several results on the best constants in different forms of maximum principles for linear elliptic equations and systems (see the monograph \cite{kresinmazbook} and the recent survey \cite{kresmazsurv}).

The case of higher order operators is quite different. Apparently, only the paper 
\cite{langermazya} by Langer and Maz’ya dealt with the question of $L^p$-dissipativity for higher order differential operators. In the case $1\leq p < \infty$, $p\neq 2$, they proved that, in the class of linear partial differential operators of order higher than two, with the domain containing $\Cspt^\infty(\Om)$, there are no generators of a contraction semigroup on $L^p(\Om)$. 
If $u$ runs over not $\Cspt^\infty(\Om)$, but only $(\Cspt^\infty(\Om))^+$ (i.e., the class of nonnegative functions of $\Cspt^\infty(\Om)$), then the result for operators with real coefficients is quite different
and really surprising:
if the operator $A$ of order $k$ has real coefficients and the integral
$$
\int_{\Om} (Au)\, u^{p-1}dx
$$
preserves the sign as $u$ runs over $(\Cspt^\infty(\Om))^+$,  then either $k=0, 1$ or $2$, or $k=4$ and $3/2\leq p\leq 3$.

Let us now give an outline of the paper. Section 
\ref{sec:functanal} presents the basic  results of Functional Analysis leading to the
concept of abstract dissipative operators. 

In Section \ref{sec:scalar} we recall our general definition of $L^p$-dissipativity of the sesquilinear form
related to a scalar second order operator. In Section 
\ref{sec:scaldiv} we give an algebraic condition we found, which provides necessary and sufficient conditions for the $L^p$-dissipativity of second order operators in divergence form, with no lower order terms.

 Section \ref{sec:pellipt} presents a review on $p$-elliptic operators, which are operators satisfying a
strengthened version of our algebraic condition.

Section \ref{sec:lower} is concerned with the $L^p$-dissipativity of operators with lower order terms. 

The topic of Section \ref{sec:elast} is the linear elasticity system. More general systems are considered
in Section \ref{sec:systems}.

In Section \ref{sec:angle} we show how the necessary and sufficient conditions we have obtained
lead to determine exactly the angle of dissipativity of certain operators. 

Section \ref{sec:max} is devoted to presents some of the results obtained by Kresin and Maz'ya 
concerning the validity of the classical maximum principle.

In Section \ref{sec:others} we briefly describe some other results we have obtained. They concern the $L^p$-dissipativity of  first order systems, of the ``complex oblique derivative'' operator and 
of a certain class of integral operators which includes the fractional powers of Laplacian
$(-\Delta)^s$, with $0<s<1$. 

Section \ref{sec:funct} discusses the concept of functional dissipativity, which we have
recently introduced. 

The final section of this paper, Section \ref{sec:conclude}, briefly shows how our 
 conditions for $L^p$-dissipativity and its strengthened variant are getting more and more important in many respects.

\section{Abstract setting}\label{sec:functanal}

Let $X$ be a (complex) Banach space. 
A \textit{semigroup of linear operators on $X$}
is a family of linear and continuous operators  
$T(t)$ ($0\leq t <\infty$) from $X$ into itself such that
$T(0)=I$,
$T(t+s)=T(t)T(s) \ (s,t\geq 0)$.

We say that $T(t)$ is a strongly continuous semigroup
 (briefly, a $C^{0}$-semigroup) if
$$
    \lim_{ t\to 0^{+}} T(t)x =x, \qquad \forall\ x\in X.
$$

The linear operator
\begin{equation}
     Ax=\lim_{t\to 0^{+}}\frac{T(t)x-x}{t}
     \label{defA}
\end{equation}
is the \textit{infinitesimal generator of the semigroup}
 $T(t)$.
 
The domain $D(A)$ of the operator $A$ (maybe not continuous) 
is the set of $x\in X$ such
that the  limit  in \eqref{defA} does exist.

If  $T(t)$ is a $C^{0}$-semigroup generated by $A$ and $u_{0}$ is a given element in  $D(A)$, 
 the function $u(t)=T(t)u_{0}$ is
solution of the evolution problem
\begin{equation}
    \begin{cases}\displaystyle \frac{du}{dt}=Au, & \text{$(t>0)$} \\
    u(0)=u_{0}\, . 
    \end{cases}
    \label{eq:pbC}
\end{equation}

 We remark the it is still possible to solve the Cauchy problem  \eqref{eq:pbC}
      when
      $u_{0}$ is an arbitrary element of $X$. In order to do that, 
      it is necessary 
      to introduce a concept of generalized solution.
      For this we refer to \cite[Ch.4]{pazy}.


One can show that if $T(t)$ is a  $C_{0}$ semigroup, then there exist two constants
    $\om\geq
     0$, $M\geq 1$ such hat
     \begin{equation}
             \Vert T(t) \Vert \leq M\, e^{\om t}, \qquad 0\leq
	     t <\infty.
         \label{eq:disC0}
     \end{equation}

 If we can choose $\om=0$ and $M=1$ in the inequality
  \eqref{eq:disC0},  we have
  $$
  \Vert T(t)\Vert \leq 1, \qquad 0\leq
	     t <\infty
  $$
  and the semigroup is said to be a contraction semigroup
  or a semigroup of contractions. If the operator $A$ is the
  generator of a $C^{0}$-semigroup of contractions, the solution of 
  the Cauchy problem \eqref{eq:pbC} satisfies the estimates
  \begin{equation}
      \Vert u(t)\Vert \leq \Vert u_{0}\Vert \qquad 0\leq
	     t <\infty .
      \label{eq:disnorm}
  \end{equation}
  
If the norm in \eqref{eq:disnorm} is the $L^\infty$ norm, we have the
classical maximum principle for parabolic equations.

The next famous Hille-Yosida Theorem characterizes the operators
which generates $C^0$  semigroups of 
contractions
\begin{theorem}
A linear operator $A$ generates a $C^0$  semigroup of 
contractions $T(t)$ if, and only if, 

(i) $A$ is closed and $D(A)$ is dense in $X$;

(ii) the resolvent set $\ro(A)$ contains $\R^+$ and the resolvent operator
$R_\la$ satisfies the inequality
$$
        \Vert R_\la\Vert \leq \frac{1}{\la}\, , 
          \quad \forall\ \la>0.
$$
\end{theorem}

Given $x\in X$, denote by $\I(x)$ the set
$$
\I(x)=\{ x^{*}\in X^{*}\ |\ \lan x^{*},x\ran =\Vert x\Vert^{2}=\Vert
x^{*}\Vert^{2}\},
$$
$X^{*}$ being the 
(topological) dual
space of $X$. 
The set $\I(x)$ is called \textit{the dual set of $x$}. 
The operator $A$ is said to be dissipative if, for any $x\in D(A)$, there exists $x^{*}\in\I(x)$ 
such that
\begin{equation}\label{eq:dissipmaincond}
  \Re\, \lan x^{*}, Ax\ran \leq 0.
\end{equation}

 Another characterization of  operators generating contractive semigroups is given
 by the equally famous Lumer-Phillips theorem:
\begin{theorem}
     If $A$ generates a  $C^0$ semigroup of contractions, then

(i) $\overline{D(A)}=X$;

(ii) $A$ is dissipative. More precisely, for any 
$x\in D(A)$, we have 
$$
\Re \lan x^*,Ax \ran \leq 0, \forall\ x^*\in \I(x);
$$

(iii) $\ro(A)\supset \R^+$.

Conversely, if 

(i') $\overline{D(A)}=X$;

(ii') $A$ is dissipative;

(iii') $\ro(A)\cap \R^+ \neq \emptyset$,

then $A$ generates a $C^0$ semigroup of contractions.
\end{theorem}

Lumer-Phillips theorem shows that in order to have a contractive
semigroup the operator $A$ has to be dissipative, i.e., inequality \eqref{eq:dissipmaincond} has to
be satisfied. If $X=L^p(\Om)$ ($1<p<\infty$) it is easy to 
see that  the dual set $\I(f)$ contains
    only the element $f^{*}$ defined by
    $$
f^{*}(x)
\begin{cases}
= \Vert f\Vert_{p}^{2-p}\overline{f(x)}\,
|f(x)|^{p-2} &\text{if $f(x)\neq 0$}\\
=0 & \text{if $f(x)=0$.}
\end{cases}
$$
and then inequality \eqref{eq:dissipmaincond} coincide
with \eqref{eq:Lpdiss}. We remark that, in the case $1<p<2$,
the integral in  \eqref{eq:Lpdiss} has to be understood
with the integrand extended by zero on the set where it vanishes.

Maz'ya and  Sobolevski\u{\i}
\cite{mazyasobolev}
obtained independently of Lumer and Phillips the same result
under the assumption that the norm of the Banach space
is G\^{a}teaux-differentiable. Their result looks as follows
\begin{theorem}
    The closed and densely defined operator $A + \la I$ has a bounded inverse for all $\la\geq 0$
    and satisfies the inequality
    $$
    \Vert [A+\la I]^{-1}\Vert \leq [\Re \la +\la_{0}]^{-1}
    $$ 
   $(\la_{0}>0)$ if and only if, for
    any $v\in D(A)$ and $f\in D(A^{*})$,
    \begin{equation*}
        \Re \lan \GA v , Av \ran \geq \la_{0}\Vert v\Vert,
    \end{equation*}
    \begin{equation*}
        \Re \lan \GA^{*}f, A^{*}f \ran \geq \la_{0}\Vert f\Vert.
        \end{equation*}
\end{theorem}
Here $\GA$ and $\GA^*$ stand for the G\^{a}teaux gradient of the norm in $B$ and in $B^*$,
respectively.
Applications to second order elliptic
operators were also given
in
\cite{mazyasobolev}.
It is interesting to note
that the paper
\cite{mazyasobolev}
was sent to the journal in 1960,
before the Lumer-Phillips paper of 1961  \cite{lumphill} appeared.

\section{Scalar second order operators with complex coefficients}
\label{sec:scalar}

In this section we describe the main results obtained in \cite{CM2005},
where we studied the $L^p$-dissipativity of scalar second order operators with complex coefficients.

To be precise we consider operators of the form
$$
Au ={\nabla^{t}}({\mathop{\mathscr A}\nolimits}\nabla u) + {\bf b}\nabla u + {\nabla^{t}}({\bf c}u)+au
$$
where ${\nabla^{t}}$ is
the divergence operator, defined
in a domain $\Om\subset \R^N$.  The coefficients
satisfy the following very general assumptions:

~---~ ${\mathop{\mathscr A}\nolimits}$ is
an $n\times n$
matrix whose entries are complex-valued measures $a^{hk}$ belonging
to $({C_{0}}({\Omega}))^{*}$. This is the dual space of ${C_{0}}({\Omega})$,
the space of
complex-valued continuous functions
with compact support contained in ${\Omega}$;

\vskip5pt
~---~ ${\bf b}=(b_{1},\ldots,b_{n})$ and ${\bf c}=(c_{1},\ldots,c_{n})$
are complex-valued vectors
with $b_{j}, c_{j}\in ({C_{0}}({\Omega}))^{*}$;

\vskip5pt
~---~ $a$ is a complex-valued scalar distribution
in $(C^{1}_{0}(\Om))^{*}$, where $C^{1}_{0}(\Om)=
C^{1}({\Omega})\cap {C_{0}}({\Omega})$.

Consider the related sesquilinear form ${\mathscr L}(u,v)$
$$
{\elle}(u,v)=\int_{{\Omega}}({\langle}{\mathop{\mathscr A}\nolimits}\nabla u,\nabla v{\rangle}-{\langle}{\bf b}\nabla u,v
\rangle
+{\langle}u,\overline{\bf c}\nabla v\rangle
-a{\langle}u,v\rangle)\,
$$
on ${C_{0}}^{1}({\Omega})\times {C_{0}}^{1}({\Omega})$.

The operator $A$ acts from $C^{1}_{0}(\Om)$ to
$(C^{1}_{0}(\Om))^{*}$ through the
relation
$$
\elle(u,v)=-\int_{{\Omega}}{\langle}Au,v\rangle
$$
for any $u,v \in C^{1}_{0}(\Om)$. The integration is
understood in the sense of distributions.

The following definition was
given in
 \cite{CM2005}.
Let $1<p<\infty$. A form ${\mathscr L}$ is
called
\textit{$L^{p}$-dissipative} if
for all $u\in C^{1}_{0}(\Om)$
\begin{equation}\label{eq:defdissorig}
\begin{gathered}
 \Re {\mathscr L}(u, |u|^{p-2}u) \geq 0,
\quad\text{if } p\geq 2, \\
\displaystyle \Re  {\mathscr L}(|u|^{p'-2}u, u) \geq 0,
\quad\text{if }
1<p< 2,
\end{gathered}
\end{equation}
where
$p'=p/(p-1)$
(we use here that $|u|^{q-2}u\in C^{1}_{0}(\Om)$ for $q\geq 2$ and
$u\in C^{1}_{0}(\Om)$).

We remark that the form ${\mathscr L}$ is $L^{p}$-dissipative if and only
if 
\begin{equation}\label{eq:elle0+}
\Re {\mathscr L}(u, |u|^{p-2}u) \geq 0
\end{equation}
for any $u\in  C_{0}^{1}({\Omega})$ such that $|u|^{p-2}u\in  C_{0}^{1}({\Omega})$.

Indeed, if $p\geq 2$,  $|u|^{p-2}u$ belongs to $C^{1}_{0}(\Om)$
for any $u\in C^{1}_{0}(\Om)$. 
If $1<p<2$, we  prove the following simple
fact:
$u\in C^{1}_{0}(\Om)$ is such that $|u|^{p-2}u$ belongs
to $C^{1}_{0}(\Om)$ if and only if 
we can write $u=\Vert v\Vert_{p'}^{2-p'} |v|^{p'-2}\overline{v}$, with $v\in C^{1}_{0}(\Om)$.

In fact, if $v$ is any function in $C^{1}_{0}(\Om)$,
then setting $u=\Vert v\Vert_{p'}^{2-p'}|v|^{p'-2}\overline{v}$, we have $u\in C^{1}_{0}(\Om)$
and $u^{*}=v$ belongs to $C^{1}_{0}(\Om)$ too.
Conversely, if $u$ is such that $|u|^{p-2}u$ belongs to $C^{1}_{0}(\Om)$,
set $v=u^{*}$. We have $v\in C^{1}_{0}(\Om)$ and $u=\Vert v\Vert_{p'}^{2-p'}|v|^{p'-2}\overline{v}$.

Therefore, if $1<p<2$, condition \eqref{eq:elle0+}
for any $u\in  C_{0}^{1}({\Omega})$ such that $|u|^{p-2}u\in  C_{0}^{1}({\Omega})$
means
$$
\Re {\mathscr L}(|v|^{p'-2}v,v)  \geq 0
$$
for any $v\in C^{1}_{0}(\Om)$.
This completes the proof of the equivalence between condition 
\eqref{eq:elle0+}
 for any $u\in  C_{0}^{1}({\Omega})$ such that $|u|^{p-2}u\in  C_{0}^{1}({\Omega})$ and  definition \eqref{eq:defdissorig}.

A first property of dissipative operators is given by the
 lemma

\begin{lemma}
If a form $\elle$ is $L^{p}$-dissipative, then
\begin{equation}\label{eq:nonnegat}
{\langle}\Re {\mathop{\mathscr A}\nolimits}\xi,\xi{\rangle}\geq 0
\quad \forall
\xi\in \R^{n}.
\end{equation}
\end{lemma}

This assertion follows from the following basic lemma wich provides 
a necessary and sufficient condition for the $L^p$-dissipativity 
of the form $\elle$ .

\begin{lemma}[\cite{CM2005}]\label{lemma:1}
A form $\elle$ is $L^{p}$-dissipative
if and only if for all $w\in C_{0}^{1}({\Om})$
\begin{align*}
{}&
\Re \int_{{\Om}}\Big[ {\langle}{\mathop{\mathscr A}\nolimits}\nabla w,\nabla w{\rangle}
-
(1-2/p){\langle}({\mathop{\mathscr A}\nolimits}-{\mathop{\mathscr A}\nolimits}^{*})\nabla(|w|),|w|^{-1}\overline{w}\nabla w{\rangle}
\\&
-
(1-2/p)^{2}{\langle}{\mathop{\mathscr A}\nolimits}\nabla(|w|),\nabla(|w|)\rangle
\Big]
+
\int_{{\Om}}{\langle}\Im  ({\bf b}+{\bf c}), \Im (\overline{w}\nabla
w)\rangle\,
\\
&
+
\int_{{\Om}} \Re ({\nabla^{t}} ({\bf b}/p - {\bf c}/p') - a
)|w|^{2}
	\geq 0.
\end{align*}
\end{lemma}

Condition \eqref{eq:nonnegat} is necessary and sufficient for the $L^2$-
dissipativity, but it is not sufficient if $p\neq 2$.

Lemma \ref{lemma:1}
implies the following sufficient condition.

\begin{corollary}[\cite{CM2005}]
Let $\alpha $ and
${\beta}$ be two real constants. If
\begin{align}
\frac{4}{p\,p'}{\langle}\Re 
{\mathop{\mathscr A}\nolimits}\xi,\xi{\rangle}
+ {\langle}\Re  {\mathop{\mathscr A}\nolimits}{\eta},{\eta}\rangle
&+2 {\langle}(p^{-1}\Im {\mathop{\mathscr A}\nolimits}+ p'^{-1}\Im {\mathop{\mathscr A}\nolimits}^{*}) \xi,{\eta}{\rangle}
\notag
\\
&
+
{\langle}\Im ({\bf b}
+ {\bf c}),{\eta}\rangle
-
2 {\langle}\Re  (\alpha {\bf b}/p - {\beta}{\bf c}/p'),\xi\rangle
\notag
\\
&+
\Re \left[{\nabla^{t}}\left((1-\alpha ){\bf b}/p - (1-{\beta}){\bf c}/p'\right) - a
\right]\geq 0
\label{polyn}
\end{align}
for any $\xi,{\eta}\in{\mathbb{R}}^{n}$,
then the form $\elle$ is $L^{p}$-dissipative.
\end{corollary}

Putting $\alpha ={\beta}=0$ in \eqref{polyn}, we find that if
\begin{eqnarray}\label{polyn2}
& \displaystyle
\frac{4}{p\,p'}{\langle}\Re 
 {\mathop{\mathscr A}\nolimits}\xi,\xi{\rangle}+ {\langle}\Re  {\mathop{\mathscr A}\nolimits}{\eta},{\eta}\rangle
+2 {\langle}(p^{-1} \Im {\mathop{\mathscr A}\nolimits}+ p'^{-1} \Im {\mathop{\mathscr A}\nolimits}^{*}) \xi,{\eta}{\rangle}
& \cr
&
+ \displaystyle {\langle} \Im ({\bf b} + {\bf c}),{\eta}\rangle
+
\Re  \left[{\nabla^{t}}\left({\bf b}/p - {\bf c}/p'\right) - a
\right]\geq 0 &
\end{eqnarray}
for any $\xi,{\eta}\in{\mathbb{R}}^{n}$,
then the form $\elle$ is $L^{p}$-dissipative.

Generally speaking,  condition \eqref{polyn2} (and the more general
condition \eqref{polyn}) is not
necessary.

\begin{example}
Let $n=2$ and
$$
{\mathop{\mathscr A}\nolimits}=\left(\begin{array}{cc}
1 & i{\gamma}\\ -i{\gamma}& 1
\end{array}\right),
$$
where ${\gamma}$ is a real constant, ${\bf b}={\bf c}=a=0$. In this case,
the polynomial \eqref{polyn2}
is given by
$$
({\eta}_{1}+{\gamma}\xi_{2})^{2}+	({\eta}_{2}-{\gamma}\xi_{1})^{2}
-({\gamma}^{2}-4/(pp'))|\xi|^{2}.
$$
For ${\gamma}^{2} >
4/(pp')$ the condition \eqref{polyn} is not satisfied,
whereas the $L^{p}$-dissipativity holds because the corresponding
operator $A$
is the Laplacian.

Note that the matrix $\Im {\mathop{\mathscr A}\nolimits}$ is not symmetric.
Below
(after Corollary \ref{cor:const}),
we give another example showing
that
the condition \eqref{polyn2} is not necessary for the $L^{p}$-dissipativity
even for symmetric
matrices $\Im {\mathop{\mathscr A}\nolimits}$.
\end{example}

\section{The operator ${\nabla^{t}}({\mathop{\mathscr A}\nolimits}\nabla u)$}\label{sec:scaldiv}

The main result proved in \cite{CM2005} concerns a scalar operator in divergence form with no lower order terms:
\begin{equation}
Au={\nabla^{t}}({\mathop{\mathscr A}\nolimits}\nabla u).
\label{eq:nolower}
\end{equation}

The following assertion
gives a necessary and sufficient
condition for the $L^{p}$-dissipativity of
the operator \eqref{eq:nolower}, where - as before -
the coefficients $a^{hk}$ belong
to
$({C_{0}}({\Omega}))^{*}$.

\begin{theorem}[\cite{CM2005}]\label{th:main}
Let
$\Im{\mathscr A}$ be symmetric, i.e.,
$\Im{\mathscr A}^{t}=$$\Im{\mathscr A}$. The form
$$
\elle(u,v)=\int_{{\Omega}}{\langle}{\mathop{\mathscr A}\nolimits}\nabla u,\nabla v\rangle
$$
is $L^{p}$-dissipative if and only if
\begin{equation}\label{eq:24}
|p-2|\, |{\langle}\Im {\mathop{\mathscr A}\nolimits}\xi,\xi\rangle| \leq 2 \sqrt{p-1}\,
{\langle}\Re {\mathop{\mathscr A}\nolimits}\xi,\xi\rangle
\end{equation}
for any $\xi\in\R^{n}$, where $|\cdot|$ denotes the total
variation.
\end{theorem}

The condition \eqref{eq:24} is understood in the
sense of comparison of measures. Of course if the coefficients 
$\{a_{hk}\}$ are complex valued $L^\infty$ functions (or more generally
$L^1_{\text{loc}}$), the
condition \eqref{eq:24} means
$$
|p-2|\, |{\langle}\Im {\mathop{\mathscr A(x)}\nolimits}\xi,\xi\rangle| \leq 2 \sqrt{p-1}\,
{\langle}\Re {\mathop{\mathscr A(x)}\nolimits}\xi,\xi\rangle
$$
for any $\xi\in\R^{n}$ and for a.e. $x\in\Om$.

When this result appeared, it was new even for operators with smooth
coefficients. For such operators it implies the contractivity of the generated
semigroup.

Note that from Theorem \ref{th:main} we immediately
derive the following well known results.

\begin{corollary}
Let $A$ be such that ${\langle}\Re 
{\mathop{\mathscr A}\nolimits}\xi, \xi {\rangle}\geq 0$ for
any $\xi\in{\mathbb{R}}^{n}$. Then

\vskip5pt
1)  $A$ is $L^{2}$-dissipative,
\vskip5pt
2) if $A$ is an operator with real coefficients, then
	$A$ is $L^{p}$-dissipative for any $p$.
\end{corollary}

The condition \eqref{eq:24} is equivalent to the positivity of
some polynomial in $\xi$ and ${\eta}$.
More exactly,
\eqref{eq:24}
is equivalent to the following condition:
\begin{equation}
\frac{4}{p\,p'}{\langle}\Re {\mathop{\mathscr A}\nolimits}\xi,\xi{\rangle}+ {\langle}\Re {\mathop{\mathscr A}\nolimits}{\eta},{\eta}\rangle
-2(1-2/p){\langle}\Im {\mathop{\mathscr A}\nolimits}\xi,{\eta}{\rangle}\geq 0
\label{eq:25}
\end{equation}
for any $\xi, {\eta}\in {\mathbb{R}}^{n}$.

More generally, if the matrix $\Im\A$ is not symmetric,
the  condition 
\begin{equation}\label{intro:form}
 \frac{4}{p\,p'}\lan \Re \A(x)\xi,\xi\ran + \lan \Re \A(x)\eta,\eta\ran
       +2 \lan (p^{-1}\Im\A(x) + p'^{-1}\Im\A^{*}(x)) \xi,\eta \ran \geq 0
\end{equation}
for almost any $x\in \Om$  and for any $\xi, \eta\in\R^{n}$  ($p'=p/(p-1)$) is only sufficient for the 
$L^p$-dissipativity.

Let us assume that either $A$ has lower order terms
or
it has no lower order terms
and $\Im {\mathop{\mathscr A}\nolimits}$ is not symmetric.
Then \eqref{eq:24} is
still necessary
for the $L^{p}$-dissipativity of $A$,
but not sufficient,
which
will be shown in
Example \ref{examp-2-2}
(cf. also
Theorem \ref{th:const} below for the case of constant
coefficients).
In other words, for such general operators the algebraic
condition \eqref{eq:25} is necessary but not sufficient,
whereas
the condition \eqref{polyn2} is sufficient, but not necessary.

\begin{example}
\label{examp-2-2}
Let $n=2$, and let ${\Omega}$ be a bounded domain. Denote by ${\sigma}$
a real function of class ${C_{0}}^{2}({\Omega})$
which does not vanish identically.
 Let
${\lambda}\in{\mathbb{R}}$.
Consider the operator \eqref{eq:nolower} with
$$
{\mathop{\mathscr A}\nolimits}=\left(\begin{array}{cc}
1 & i{\lambda}{\partial}_{1}({\sigma}^{2}) \\ -i{\lambda}{\partial}_{1}({\sigma}^{2}) & 1
\end{array}\right),
$$
i.e.,
$$
Au={\partial}_{1}({\partial}_{1}u+i{\lambda}{\partial}_{1}({\sigma}^{2})\, {\partial}_{2}u) +
{\partial}_{2}(-i{\lambda}{\partial}_{1}({\sigma}^{2})\,
{\partial}_{1}u+{\partial}_{2}u),
$$
where ${\partial}_{i}={\partial}/{\partial}x_{i}$ ($i=1,2$).
By definition, we have $L^{2}$-dissipativity if and only if
$$
\Re \int_{{\Omega}}(({\partial}_{1}u+i{\lambda}{\partial}_{1}({\sigma}^{2})\, {\partial}_{2}u)
{\partial}_{1}\overline{u}
+ (-i{\lambda}{\partial}_{1}({\sigma}^{2})\,
{\partial}_{1}u+{\partial}_{2}u){\partial}_{2}\overline{u})\, dx \geq 0
$$
for any $u\in{C_{0}}^{1}({\Omega})$, i.e., if and only if
$$
\int_{{\Omega}}|\nabla u|^{2}dx -2{\lambda}\int_{{\Omega}}{\partial}_{1}({\sigma}^{2})
\Im ({\partial}_{1}\overline{u}\,{\partial}_{2}u)\, dx \geq 0
$$
for any $u\in{C_{0}}^{1}({\Omega})$.
Taking $u={\sigma}\, \exp(itx_{2})$ ($t\in{\mathbb{R}}$), we obtain, in particular,
\begin{equation}
t^{2}\int_{{\Omega}}{\sigma}^{2}dx -t
\lambda
\int_{{\Omega}}({\partial}_{1}({\sigma}^{2}))^{2} dx + \int_{{\Omega}}|\nabla {\sigma}|^{2}dx \geq 0.
\label{eq:inpartic}
\end{equation}
Since
$$
\int_{{\Omega}}({\partial}_{1}({\sigma}^{2}))^{2} dx > 0,
$$
we can choose ${\lambda}\in{\mathbb{R}}$ so that
\eqref{eq:inpartic} is impossible
for all $t\in{\mathbb{R}}$.
Thus, $A$ is not $L^{2}$-dissipative, although
\eqref{eq:24} is satisfied.
Since $A$ can be written as
$$
Au=\Delta u - i{\lambda}({\partial}_{21}({\sigma}^{2})\,
{\partial}_{1}u - {\partial}_{11}({\sigma}^{2})\, {\partial}_{2}u),
$$
this example shows
that \eqref{eq:24} is not sufficient for the $L^{2}$-dissipativity
of an operator with
lower order terms,
even if $\Im {\mathop{\mathscr A}\nolimits}$ is symmetric.
\end{example}


\section{The $p$-ellipticity}\label{sec:pellipt}

 Let us consider the class of operators  
 \begin{equation}\label{eq:lowerterms}
 Au= \nabla(\A \nabla u)
+  { b}\nabla u + \nabla(c u)+au.
\end{equation}
 with $L^\infty$ coefficients, such that
the form \eqref{intro:form} is not merely
non-negative, but strictly positive, i.e., there exists $\kappa>0$ such that
\begin{equation}\label{intro:formk}
\begin{gathered}
 \frac{4}{p\,p'}\lan \Re \A(x)\xi,\xi\ran + \lan \Re \A(x)\eta,\eta\ran
       +2 \lan (p^{-1}\Im\A(x) + p'^{-1}\Im\A^{*}(x)) \xi,\eta \ran \\
        \geq \kappa(|\xi|^{2}+|\eta|^{2})
       \end{gathered}
\end{equation}
for almost any $x\in \Om$  and for any $\xi, \eta\in\R^{n}$. 
The class of operators \eqref{eq:lowerterms} whose principal part satisfies \eqref{intro:formk} and which  could be called (\textit{strongly}) \textit{$p$-elliptic}, 
was recently considered by several authors.

Carbonaro and Dragi\v{c}evi\'c \cite{carbonaro,carbonaro2} 
 showed the validity of a so called (dimension free) bilinear embedding. Their main result is the following
 \begin{theorem}[\cite{carbonaro}]
   Let $P_{t}^{A}=\exp \left(-t L_{A}\right), t>0$ and
let $p>1$. Suppose that the matrices $A, B$ are  $p$-elliptic. Then for all $f, g \in \Cspt^{\infty}\left(\mathbb{R}^{n}\right)$ we have
\begin{equation}\label{eq:(1.12)}
\int_{0}^{\infty} \int_{\R^{n}}\left|\nabla P_{t}^{A} f(x)\right|\left|\nabla P_{t}^{B} g(x)\right| d x d t \leq C \|f\|_{p}\|g\|_{p'}
\end{equation}
with constant depending on ellipticity parameters, but not dimension.  
\end{theorem}

We note that if $A$ and $B$ are real accretive matrices then
 \eqref{eq:(1.12)} holds for the full range of exponents $p\in (1,\infty)$.

In a series of papers  \cite{DP20191,DP20192,DP20201,DP20202}
 Dindo\v{s} and Pipher
proved several results concerning the $L^p$ solvability of the Dirichlet problem. Their 
result hinges on a regularity property for the  solutions of
the Dirichlet problem for the equation
\begin{equation}\label{eq:dind}
\partial_{i}\left(a_{i j}(x) \partial_{j} u\right)+b_{i}(x) \partial_{i} u=0\, ,
\end{equation}
given by the next result 

\begin{lemma}[\cite{DP20191}, p.269]
Let the matrix $A$ be $p$-elliptic for $p \geq 2$ and let $B$ have 
coefficients  $B_{i} \in L_{\text {\rm loc }}^{\infty}(\Omega)$ satisfying the condition
\begin{equation}\label{eq:(1.9)}
\left|B_{i}(x)\right| \leq K(\delta(x))^{-1}, \quad \forall x \in \Omega\, ,
\end{equation}
where the constant $K$ is uniform, and $\delta(x)$ denotes the distance of $x$ to the boundary of $\Omega$.  Suppose that $u \in W_{\text{loc}}^{1,2}(\Omega)$ is a weak solution 
  of the equation \eqref{eq:dind} in $\Omega$, an open subset of $\R^{n}$. Then, for any ball $B_{r}(x)$ with $r<\delta(x) / 4$,
$$
\int_{B_{r}(x)}|\nabla u(y)|^{2}|u(y)|^{p-2} d y \leq C_1 r^{-2} \int_{B_{2 r}(x)}|u(y)|^{p} d y
$$
and
$$
\left(\int_{B_{r}(x)}|u(y)|^{q} d y\right)^{1 / q} \leq C_2 \left(\int_{B_{2 r}(x)}|u(y)|^{2} d y\right)^{1 / 2}
$$
for all $q \in\left(2, \frac{n p}{n-2}\right]$ when $n>2$, and where
$C_1, C_2$ depend only  on $p$-ellipticity constants and $K$ of \eqref{eq:(1.9)}. When $n=2, q$ can be any number in $(2, \infty)$. In particular, $|u|^{(p-2) / 2} u$ belongs to $W_{\text \rm {loc }}^{1,2}(\Omega)$.
\end{lemma}

Dindo\v{s} and Pipher used this result in an iterative procedure, which is similar to Moser’s iteration scheme (used in Moser's proof  of
the celebrated De Giorgi–Nash–Moser regularity theorem for real divergence form elliptic equations). 
Differently from the real coefficients case, where the procedure can be
applied for any $p$ and leads to the boundedness of the solution (and then to its H\"olderianity),  here the iteration scheme can be applied only 
up to a threshold determined by the $p$-ellipticity of the operator.  
This is sufficient to obtain an higher integrability of the solution. 

Dindo\v{s} and Pipher uses  this regularity result in the study of the  existence 
for the Dirichlet problem
\begin{equation}\label{eq:D(4.1)}
\begin{cases}\partial_{i}\left(a_{i j}(x) \partial_{j} u\right)+b_{i}(x) \partial_{i} u =0 & \text {in } \Omega \\ u(x)=f(x) &  \text {a.e. on }  \partial \Omega \\ \widetilde{N}_{2, a}(u) \in L^{p}(\partial \Omega) & \end{cases}
\end{equation}
where $f$ is in $L^{p}(\partial \Omega)$. Here $a>0$ is a fixed  parameter and $\widetilde{N}_{2, a}(u)$ is a nontangential maximal function
 defined using $L^p$ averages over balls 
 $$
 \widetilde{N}_{2, a}(u)(y) = \sup_{x\in \Gamma_{a}(y)} 
\left(\dashint_{B_{\delta(x) / 2}(x)}|u(z)|^{2} d z\right)^{1 / 2}
 $$
($y\in\partial\Omega$) where the barred integral indicates the average and 
$\Gamma_{a}(y)$  
is a cone of aperture $a$.
To be precise, 
they say that the Dirichlet problem \eqref{eq:D(4.1)} is solvable for a given 
$p \in(1, \infty)$ if there exists a $C=C(p, \Omega)>0$ such that for all complex valued boundary data $f \in L^{p}(\partial \Omega) \cap B_{1 / 2}^{2,2}(\partial \Omega)$ the unique ``energy solution'' satisfies the estimate
$$
\left\|\widetilde{N}_{2, a}(u)\right\|_{L^{p}(\partial \Omega)} \leq C\|f\|_{L^{p}(\partial \Omega)}\, .
$$
Since the space $\dot{B}_{1 / 2}^{2,2}(\partial \Omega) \cap L^{p}(\partial \Omega)$ is dense in $L^{p}(\partial \Omega)$ for each $p \in(1, \infty)$,  there exists a unique continuous extension of the solution operator $f \mapsto u$ to the whole space $L^{p}(\partial \Omega)$, with $u$ such that $\widetilde{N}_{2, a}(u) \in L^{p}(\partial \Omega)$ and the accompanying estimate $\left\|\widetilde{N}_{2, a}(u)\right\|_{L^{p}(\partial \Omega)} \leq C\|f\|_{L^{p}(\partial \Omega)}$ is valid.

Their results have been extended by Dindo\v{s}, Li and Pipher to systems
and in particular to elasticity  in \cite{DLP}.

We mention  that  - as Carbonaro and Dragi\v{c}evi\'c \cite{carbonaro} show -
the $p$-ellipticity comes into play also in the study of the
convexity of power functions (Bellman functions) and
in the holomorphic functional calculus. 

Egert \cite{egert} shows that the $p$-ellipticity condition
 implies extrapolation to a holomorphic semigroup on 
 Lebesgue spaces in a $p$-dependent range of exponents.

Finally we remark that, if the partial differential operator has no lower order terms, the concepts of $p$-ellipticity
and strict  $L^p$-dissipativity coincide. By strict  $L^p$-dissipativity 
we mean that there exists $\kappa>0$ such that
$$
\Re \int_{\Om} \lan \A^{hk} \de_{k} u, \de_{h}(|u|^{p-2} u)\ran\, dx \geq 
\kappa \int_{\Om} | \nabla(|u|^{(p-2)/2} u)|^2 dx
$$
 for any $u\in  C_{0}^{1}({\Omega})$ such that $|u|^{p-2}u\in  C_{0}^{1}({\Omega})$.

\label{rem:dissconk}
It is worthwhile to remark that, if the partial differential operator has no lower order terms, the concepts of $p$-ellipticity
and strict  $L^p$-dissipativity coincide. 
 One can prove that the operator $A$ is strict  $L^p$-dissipative, i.e., $p$-elliptic,  if and only if there exists $\kappa>0$ such that
 $A-\kappa\Delta$ is $L^p$-dissipative.

\section{$L^p$-dissipativity for operators with lower order terms}
\label{sec:lower}

Generally speaking, it is impossible
to obtain an algebraic
characterization for
an
operator with
lower order terms.
Indeed, let us consider,
for example, the operator
$$
Au=\Delta u + a(x)u
$$
in a bounded domain ${\Omega}\subset{\mathbb{R}}^{n}$
with zero Dirichlet
boundary data. Denote by ${\lambda}_{1}$
the first eigenvalue of the Dirichlet problem for
the Laplace equation
in ${\Omega}$. A sufficient condition for
the $L^{2}$-dissipativity of $A$
has the form
$\Re  a\leq {\lambda}_{1}$, and we cannot
give an algebraic characterization of
${\lambda}_{1}$.

Consider, as another example, the operator
\begin{equation}\label{eq:Deltamu}
A =\Delta + \mu
\end{equation}
where $\mu$ is a nonnegative Radon measure on ${\Omega}$.
The operator $A$ is $L^{p}$-dissipative if and only if
\begin{equation}
\int_{{\Omega}}|w|^{2}d\mu \leq \frac{4}{pp'}\int_{{\Omega}}|\nabla w|^{2}dx
\label{eq:schro}
\end{equation}
for any $w\in C_{0}^{\infty}({\Omega})$ (cf. Lemma
\ref{lemma:1}).
Maz'ya \cite{mazya62, mazya64, mazyasimon}
proved
that the following condition is sufficient for
\eqref{eq:schro}:
\begin{equation}
\frac{\mu(F)}{\hbox{\rm cap}_{{\Omega}}(F)} \leq \frac{1}{pp'}
\label{eq:schro1}
\end{equation}
for all compact set $F\subset {\Omega}$
and the following condition is necessary:
\begin{equation}
\frac{\mu(F)}{\hbox{\rm cap}_{{\Omega}}(F)} \leq \frac{4}{pp'}
\label{eq:schro2}
\end{equation}
for all compact set $F\subset {\Omega}$.
Here, $\hbox{\rm cap}_{{\Omega}}(F)$ is the capacity of $F$
relative to
${\Omega}$, i.e.,
$$
\hbox{\rm cap}_{{\Omega}}(F)=\inf\Big\{\int_{{\Omega}}|\nabla u|^{2}dx\, :\
u\in {C_{0}}^{\infty}({\Omega}),\ u\geq 1\ \hbox{\rm on}\ F\Big\}.
$$
The condition \eqref{eq:schro1}
is not necessary and the condition \eqref{eq:schro2} is not sufficient.

It must be pointed out that a theorem by Jaye, Maz'ya and Verbitsky
can provide a necessary and sufficient condition of a different nature for the $L^p$-dissipativity
of operator \eqref{eq:Deltamu}. In fact  in \cite{JMV} they proved the following result
\begin{theorem}
  Let $\Omega$ be an open set, and let $\sigma \in  (\Cspt^\infty(\Om))'$ be a real valued distribution. In addition, let $\mathcal{A}$ be a symmetric matrix function defined on $\Omega$ satisfying the conditions
$$
\text { (1.4) } \quad m|\xi|^{2} \leq \mathcal{A}(x) \xi \cdot \xi, \quad \text { and } \quad|\mathcal{A}(x) \xi| \leq M|\xi|, \quad \text { for all } \xi \in \R^{n} \setminus\{0\}.
  $$
  Then 
$$
\left\langle\sigma, h^{2}\right\rangle \leq \int_{\Omega}(\mathcal{A} \nabla h) \cdot \nabla h\,  d x
$$
 holds for all $h \in \Cspt^{\infty}(\Omega)$ if and only if there exists a vector field $\vec{\Gamma} \in L_{\text{loc}}^{2}(\Omega)$ so that:
$$
\sigma \leq \dive(\mathcal{A} \vec{\Gamma})-(\mathcal{A} \vec{\Gamma}) \cdot \vec{\Gamma} \quad \text { in }(\Cspt^\infty(\Om))' .
$$  
\end{theorem}

Keeping in mind that the operator \eqref{eq:Deltamu} is $L^p$-dissipative if and only if 
\eqref{eq:schro} holds for any $w\in \Cspt^{\infty}(\Omega)$, 
by taking $\si=\mu$, $\mathcal{A} =(4/(pp'))\mathcal{I}$, we find immediately that 
\eqref{eq:Deltamu} is $L^p$-dissipative  if and only if there exists a vector field
 $\vec{\Gamma} \in L_{\text{loc}}^{2}(\Omega)$ such that
 $$
 \mu \leq \frac{4}{p\, p'} \left( \dive \vec{\GA} - |\vec{\GA}|^2 \right)
 $$
in the sense of distributions.

In the case of an operator with constant coefficients  and lower order terms, 
 we have found a  necessary and
sufficient condition.
Consider the operator
\begin{equation}\label{eq:genform}
Au=\nabla^{t}({\mathop{\mathscr A}\nolimits}\nabla u) + {\bf b}\nabla u + au
\end{equation}
with constant complex coefficients. Without loss of generality,
we
can assume that the matrix ${\mathop{\mathscr A}\nolimits}$ is symmetric.

The following assertion
provides
a necessary and sufficient condition
for the $L^{p}$-dissipativity of the operator $A$.

\begin{theorem}[\cite{CM2005}]
\label{th:const}
Suppose that
${\Omega}$ is an open set in ${\mathbb{R}}^{n}$
which contains balls of
arbitrarily large radius. The operator \eqref{eq:genform}
is $L^{p}$-dissipative if and only if there exists a real constant
vector $V$ such that
\begin{align*}
&2\Re {\mathop{\mathscr A}\nolimits}V+\Im  {\bf b}=0,
\\
&\Re  a + {\langle}\Re {\mathop{\mathscr A}\nolimits}V,V{\rangle}\leq 0
\end{align*}
and for  any
$\xi\in{\mathbb{R}}^{n}$
\begin{equation}\label{eq:V1}
|p-2|\, |{\langle}\Im {\mathop{\mathscr A}\nolimits}\xi,\xi\rangle| \leq 2 \sqrt{p-1}\,
{\langle}\Re {\mathop{\mathscr A}\nolimits}\xi,\xi\rangle.
\end{equation}
\end{theorem}

If the matrix $\Re {\mathop{\mathscr A}\nolimits}$ is not
singular, the following assertion holds.

\begin{corollary}[\cite{CM2005}; cf. also
\cite{KM1994}]
\label{cor:const}
Suppose that
${\Omega}$ is
an open set in ${\mathbb{R}}^{n}$ which contains balls of
arbitrarily large radius.
Assume that
the matrix $\Re {\mathop{\mathscr A}\nolimits}$
is
not singular. The operator $A$
is $L^{p}$-dissipative if and only if \eqref{eq:V1} holds and
\begin{equation}
4\Re  a \leq -{\langle}(\Re {\mathop{\mathscr A}\nolimits})^{-1}\Im  {\bf b},
\Im  {\bf b}\rangle.
\label{eq:newV2}
\end{equation}
\end{corollary}

Now, we can show that the condition
\eqref{polyn} is not necessary for the
$L^{p}$-dissipativity, even if
the matrix $\Im {\mathop{\mathscr A}\nolimits}$ is symmetric.

\begin{example}
Let $n=1$,  and let $\Omega ={\mathbb{R}}^{1}$. Consider the operator
$$
\left( 1 + 2\, \frac{\sqrt{p-1}}{p-2}\, i\right) u'' +2i u' -u,
$$
where $p\neq 2$ is fixed. The conditions
\eqref{eq:V1} and \eqref{eq:newV2} are satisfied,
and this operator is $L^{p}$-dissipative
in view of Corollary \ref{cor:const}.

On the other hand, the polynomial in \eqref{polyn2}
has the form
$$
\left(2\, \frac{\sqrt{p-1}}{p}\, \xi -{\eta}\right)^{2}
+2{\eta}+1,
$$
i.e.,
it is not nonnegative for any $\xi,{\eta}\in{\mathbb{R}}$.
\end{example}

Recently Maz’ya and Verbitsky \cite{MV2019} (see also
\cite{MV2020}) gave necessary and sufficient conditions for the accretivity of a second order partial differential operator $\mathcal{L}$ containing lower order terms, in the case of Dirichlet data. We observe that the accretivity of $\mathcal{L}$ is equivalent to the $L^2$-dissipativity of $-\mathcal{L}$.

Their result concern second order operators with distributional coefficients
\begin{equation}\label{eq:MVelle}
\mathcal{L}=\operatorname{div}(A \nabla \cdot)+{\bf b} \cdot \nabla+c
\end{equation}
 where $A \in ((\Cspt^\infty(\Om))')^{n \times n}$, ${\bf b} \in ((\Cspt^\infty(\Om))')^{n}$ and $c \in (\Cspt^\infty(\Om))'$ are complex-valued.

Given $A=\left\{a_{j k}\right\} \in ((\Cspt^\infty(\Om))')^{n \times n}$, we 
denote by $A^{s}$ and $A^{c}$  its symmetric  part and skew-symmetric part  respectively. The accretivity property for $-\mathcal{L}$ can be characterized in terms of the following real-valued expressions:
\begin{equation}\label{eq:(4.1)}
P=\Re A^{s}, \quad {\bf d}=\frac{1}{2}\left[\Im {{\bf b}}-\operatorname{Div}\left(\Im   A^{c}\right)\right], \quad \sigma=\Re  c-\frac{1}{2} \operatorname{div}(\Re  {\bf b})\, .
\end{equation}
We note that $P=\left\{p_{i k}\right\} \in ((\Cspt^\infty(\Om))')^{n \times n}, {\bf d}=\left\{d_{j}\right\} \in ((\Cspt^\infty(\Om))')^{n}$, and $\sigma \in (\Cspt^\infty(\Om))'$. 

Moreover, in order that $-\mathcal{L}$ be accretive, the matrix $P$ must be nonnegative definite, i.e., $P \xi \cdot \xi \geq 0$ in $(\Cspt^\infty(\Om))'$ for all $\xi \in \mathbb{R}^{n} .$ In particular, each $p_{j j}$ $(j=1, \ldots, n)$ is a nonnegative Radon measure.

The characterization of accretive operators $-\mathcal{L}$ is given in the following criterion obtained in \cite[Proposition 2.1]{MV2019}

\begin{theorem}
    Let $\mathcal{L}$ be the operator \eqref{eq:MVelle}. Suppose that $P, {\bf d}$, and $\sigma$ are defined by \eqref{eq:(4.1)}.
The operator $-\mathcal{L}$ is accretive if and only if $P$ is a nonnegative definite matrix, and the following two conditions hold:
$$
[h]_{\mathcal{H}}^{2}=\langle P \nabla h, \nabla h\rangle-\langle\sigma h, h\rangle \geq 0
$$
for all real-valued $h \in \Cspt^{\infty}(\Omega)$, and the commutator inequality
\begin{equation}\label{eq:comm1}
|\langle {\bf d}, u \nabla v-v \nabla u\rangle| \leq[u]_{\mathcal{H}}[v]_{\mathcal{H}}
\end{equation}
holds for all real-valued $u, v \in \Cspt^{\infty}(\Omega)$.
\end{theorem}

Under some mild restrictions on $\mathcal{H}$, the ``norms''
 $[u]_{\mathcal{H}}$ and $[v]_{\mathcal{H}}$ on the right-hand side of \eqref{eq:comm1} can be replaced, up to a constant multiple, with the corresponding Dirichlet norms
 \begin{equation}\label{eq:comm2}
|\langle {\bf d}, u \nabla v-v \nabla u\rangle| \leq C\|\nabla u\|_{L^{2}(\Omega)}\|\nabla v\|_{L^{2}(\Omega)}
\end{equation}
where $C>0$ is a constant which does not depend on real-valued $u, v \in \Cspt^{\infty}(\Omega)$. 
This leads to explicit criteria of accretivity (see \cite[Section 4]{MV2020} for the details). Indeed Maz'ya and Verbitsky have found necessary and sufficient conditions for the validity of commutator inequality \eqref{eq:comm2}. For example, when $\Om=\R^n$ and ${\bf d}$ has $L^1_{\text{loc}}$ 
components, they prove the following result

\begin{theorem}
    Let ${\bf d} \in [L_{\text {loc}}^{1}(\mathbb{R}^{n})]^n, n \geq 2$. The inequality 
    \begin{equation}\label{eq:comm3}
\left| \int_{\R^n}\langle {\bf d}, u \nabla v-v \nabla u\rangle  dx \right| \leq C\|\nabla u\|_{L^{2}(\R^n)}\|\nabla v\|_{L^{2}(\R^n)}
\end{equation}
holds for any real-valued $u, v \in \Cspt^{\infty}(\R^n)$if and only if
\begin{equation}\label{eq:1.10}
{\bf d}={\bf c}+\operatorname{Div} F
\end{equation}
where $F \in \mathrm{BMO}\left(\mathbb{R}^{n}\right)^{n \times n}$ is a skew-symmetric matrix field, and ${\bf c}$ satisfies the condition
\begin{equation}\label{eq:1.11}
\int_{\mathbb{R}^{n}}|{\bf c}|^{2}|u|^{2} d x \leq C\|\nabla u\|_{L^{2}\left(\R^{n}\right)}^{2}
\end{equation}
where the constant $C$ does not depend on $u \in \Cspt^{\infty}\left(\mathbb{R}^{n}\right)$.
Moreover, if \eqref{eq:comm3} holds, then \eqref{eq:1.10} is valid with ${\bf c}=\nabla \Delta^{-1} 
(\dive \left.{\bf d}\right)$ satisfying \eqref{eq:1.11}, and $F=\Delta^{-1}(\operatorname{Curl} {\bf d}) \in \operatorname{BMO}\left(\mathbb{R}^{n}\right)^{n \times n} .$

In the case $n=2$, necessarily ${\bf c}=0$, and ${\bf d}=\left(-\partial_{2} f, \partial_{1} f\right)$ with $f \in$ $\operatorname{BMO}\left(\mathbb{R}^{2}\right)$ in the above statements.
\end{theorem}

\section{Elasticity}\label{sec:elast}

Consider the classical operator of
two-dimensional elasticity
\begin{equation}
Eu=\Delta u + (1-2{\nu})^{-1}\nabla {\nabla^{t}} u,
\label{opelast}
\end{equation}
where ${\nu}$ is the Poisson ratio.
As is known,
$E$
is strongly elliptic if and only if either ${\nu}>1$ or ${\nu}<1/2$.
To obtain a necessary and sufficient condition
for the $L^{p}$-dissipativity
of this  elasticity system,
we formulate some facts about
systems of partial differential equations
of the form
\begin{equation}
A={\partial}_{h} ({\mathop{\mathscr A}\nolimits}^{hk}(x){\partial}_{k}),
\label{eq:A2}
\end{equation}
where ${\mathop{\mathscr A}\nolimits}^{hk}(x)=\{a^{hk}_{ij}(x)\}$ are
$m\times m$ matrices
whose entries
are complex locally integrable functions defined in
an arbitrary domain ${\Omega}$ of ${\mathbb{R}}^{n}$
$(1\leq i,j\leq
m,\ 1\leq h,k \leq n)$.

\begin{lemma}[\cite{CM2006}]
\label{lemma:5}
An operator $A$ of the form
\eqref{eq:A2} is $L^{p}$-dissipative in
${\Omega}\subset{\mathbb{R}}^{n}$ if and only if
\begin{align*}
\int_{{\Omega}}
\Big(\operatorname{Re} {\langle}
&{\mathop{\mathscr A}\nolimits}^{hk}
{\partial}_{k}w,
{\partial}_{h}w\rangle
\\
&
-(1-2/p)^{2}|w|^{-4}\operatorname{Re}{\langle}{\mathop{\mathscr A}\nolimits}^{hk} w,w{\rangle}\operatorname{Re}
\langle
w,{\partial}_{k}w{\rangle}\operatorname{Re} {\langle}w,{\partial}_{h}w{\rangle}
\\[2pt]
&
-(1-2/p)|w|^{-2}\operatorname{Re}\big({\langle}{\mathop{\mathscr A}\nolimits}^{hk} w,{\partial}_{h}w{\rangle}\operatorname{Re} {\langle}w,{\partial}_{k}w\rangle
\\
&-{\langle}{\mathop{\mathscr A}\nolimits}^{hk} {\partial}_{k}w,w{\rangle}\operatorname{Re} {\langle}w,
{\partial}_{h}w\rangle\big)
\Big) dx
\geq 0
\end{align*}
for any $w\in ({C_{0}}^{1}({\Omega}))^{m}$.
\end{lemma}

In the
case $n=2$,
Lemma \ref{lemma:5}
yields a necessary algebraic condition.

\begin{theorem}[\cite{CM2006}]
\label{th:5}
Let ${\Omega}$ be a domain of ${\mathbb{R}}^{2}$.
If an operator $A$ of the form \eqref{eq:A2} is $L^{p}$-dissipative,
then
\begin{align*}
{}&
\operatorname{Re} {\langle}({\mathop{\mathscr A}\nolimits}^{hk}(x)\xi_{h}\xi_{k}){\lambda},{\lambda}{\rangle} -
(1-2/p)^{2}\operatorname{Re}{\langle}({\mathop{\mathscr A}\nolimits}^{hk}(x)\xi_{h}\xi_{k}){\omega},{\omega}{\rangle}(\operatorname{Re}
\langle
{\lambda},{\omega}\rangle)^{2}
\\[2pt]
&
-(1-2/p)\operatorname{Re}({\langle}({\mathop{\mathscr A}\nolimits}^{hk}(x)\xi_{h}\xi_{k}){\omega},{\lambda}\rangle
	 -{\langle}({\mathop{\mathscr A}\nolimits}^{hk}(x)\xi_{h}\xi_{k}){\lambda},{\omega}\rangle)
	 \operatorname{Re} {\langle}{\lambda},{\omega}{\rangle}
\geq 0
\end{align*}
for almost every $x\in{\Omega}$ and
for any $\xi\in{\mathbb{R}}^{2}$, ${\lambda},\, {\omega}\in {\mathbb{C}}^{m}$, $|{\omega}|=1$.
\end{theorem}

Based on
Lemma \ref{lemma:5} and Theorem \ref{th:5},  it is
possible to obtain
the following criterion
for the $L^{p}$-dissipativity of the two-dimensional elasticity system.

\begin{theorem}[\cite{CM2006}]
\label{th:47-52}
The operator
\eqref{opelast} is $L^{p}$-dissipative if
and only if
\begin{equation}\label{condvecchiaintro}
\left(\frac{1}{2}-\frac{1}{p}\right)^{2} \leq
\frac{2({\nu}-1)(2{\nu}-1)}
	{(3-4{\nu})^{2}}.
\end{equation}
\end{theorem}

By Theorems \ref{th:5} and \ref{th:47-52},
it is easy to compare
$E$ and $\Delta$ from the point of view of $L^{p}$-dissipativity.

\begin{corollary}[\cite{CM2006}]
There exists $k>0$ such that $E-k\Delta$ is $L^{p}$-dissipative if
and only if
$$
\left(\frac{1}{2}-\frac{1}{p}\right)^{2} <
\frac{2({\nu}-1)(2{\nu}-1)}
		 {(3-4{\nu})^{2}}\, .
$$

There exists $k<2$ such that $k\Delta-E$ is $L^{p}$-dissipative
if and only if
$$
	 \left(\frac{1}{2}-\frac{1}{p}\right)^{2} <
	 \frac{2{\nu}(2{\nu}-1)}
	 {(1-4{\nu})^{2}}
	 \, .
$$
\end{corollary}

As remarkd at p.\pageref{rem:dissconk} for scalar operators, this 
is equivalent to say that $E$ is strict  $L^p$-dis\-si\-pa\-tive, i.e., $E$ is $p$-elliptic.  The last result was
recently extended to variable Lam\'e parameters by Dindo\v{s}, Li and Pipher \cite{DLP}.
It must be pointed out that  these authors introduce an auxiliary function  $r(x)$  (see \cite[pp.390--391]{DLP}) which
generates some first order terms in the partial differential operator. In the definition of $p$-ellipticity these terms do not play any role,
while they have some role in the  dissipativity. Therefore our and their results do not seem to be completely
equivalent.

In  \cite{CM2013}
 we showed that condition \eqref{condvecchiaintro}
is necessary for the $L^{p}$-dissipativity of operator
\eqref{opelast} in any dimension, even when the Poisson ratio is not constant.
At the present it is not known if condition 
\eqref{condvecchiaintro} is  also sufficient
for the $L^{p}$-dissipativity of elasticity operator
for $n> 2$, in particular for $n=3$ (see \cite[Problem 43]{mazyaunsolved}). 
Nevertheless, in the same paper, we gave a more strict 
explicit condition which is sufficient 
for the $L^{p}$-dissipativity of \eqref{opelast}.
Indeed we proved that if
 $$
        (1-2/p)^{2}\leq
\begin{cases}
\displaystyle\frac{1-2\n}{2(1-\n)} & \text{if}\ \n<1/2\\
\\
\displaystyle\frac{2(1-\n)}{1-2\n}& \text{if}\ \n > 1,
\end{cases}
$$
 then  the operator \eqref{opelast} is $L^p$-dissipative.

In \cite{CM2013} we gave   also necessary and
sufficient conditions for a weighted $L^p$-negativity of the
Dirichlet-Lam\'e operator, i.e.
for the validity of the inequality
\begin{equation}\label{eq:introneg}
    \int_{\Om}(\Delta u + (1-2\nu)^{-1}\nabla \dive u)\,
    |u|^{p-2}u\, \frac{dx}{|x|^\alpha} \leq 0
    \end{equation}
under the condition that the vector $u$ is rotationally invariant, i.e. $u$ depends only
on $\ro=|x|$ and $u_\ro$ is the only nonzero spherical
component of $u$. Namely we showed that \eqref{eq:introneg} holds for any such $u$
belonging to $(\Cspt^{\infty}(\R^{N}\setminus\{0\}))^{N}$ if and only if
$$
  -(p-1)(n+p'-2) \leq \alpha \leq n+p-2.
$$

\section{A Class of Systems of Partial Differential
Operators}\label{sec:systems}

In this section, we consider
systems of partial differential operators of the form
\begin{equation}
Au = {\partial}_{h}({\mathop{\mathscr A}\nolimits}^{h}(x){\partial}_{h}u),
\label{defspde}
\end{equation}
where ${\mathop{\mathscr A}\nolimits}^{h}(x)=\{a_{ij}^{h}(x)\}$ ($i,j=1,\ldots,m$)
are matrices
with complex locally integrable entries defined in a domain
${\Omega}\subset{\mathbb{R}}^{n}$ ($h=1,\ldots,n$).
Note that the elasticity system is not a system of this kind.

To characterize the $L^{p}$-dissipativity of such
operators, one can reduce the consideration
to the one-dimensional case.
Auxiliary facts are given in the following two subsections.

Langer and Maz'ya considered the $L^p$-dissipativity of weakly coupled systems
in \cite{langer}.

\subsection{Dissipativity of
systems of ordinary differential
equations}

In this subsection, we consider the operator
\begin{equation}
A u = ({\mathop{\mathscr A}\nolimits}(x)u')',
\label{defAord}
\end{equation}
where ${\mathop{\mathscr A}\nolimits}(x)=\{a_{ij}(x)\}$ ($i,j=1,\ldots,m$) is a matrix with
complex locally integrable entries defined in a bounded or
unbounded
interval $(a,b)$.
The corresponding sesquilinear form $\elle(u,w)$
takes the form
$$
\elle(u,w)=\int_{a}^{b}{\langle}{\mathop{\mathscr A}\nolimits}u',w'{\rangle}\, dx .
$$

\begin{theorem}[\cite{CM2006}]
\label{th:0}
The operator $A$ is $L^{p}$-dissipative if and only if
\begin{align*}
{}&
\Re  {\langle}{\mathop{\mathscr A}\nolimits}(x) {\lambda},{\lambda}\rangle-(1-2/p)^{2}\Re  {\langle}{\mathop{\mathscr A}\nolimits}(x){\omega},
	 {\omega}{\rangle}(\Re  \langle{\lambda},{\omega}\rangle)^{2}
	 \\[2pt]
&-
(1-2/p)\Re ({\langle}{\mathop{\mathscr A}\nolimits}(x){\omega},{\lambda}{\rangle}-{\langle}{\mathop{\mathscr A}\nolimits}(x){\lambda},{\omega}\rangle)
	 \Re  {\langle}{\lambda},{\omega}{\rangle} \geq 0
\end{align*}
for almost every $x\in(a,b)$ and for any
${\lambda},{\omega}\in{\mathbb{C}}^{m}$, $|{\omega}|=1$.
\end{theorem}

This theorem implies the following assertion.

\begin{corollary}[\cite{CM2006}]
If the operator $A$ is $L^{p}$-dissipative, then
$$
\Re  {\langle}{\mathop{\mathscr A}\nolimits}(x){\lambda},{\lambda}{\rangle}\geq 0
$$
for almost every $x\in(a,b)$ and for any ${\lambda}\in{\mathbb{C}}^{m}$.
\end{corollary}

As a consequence of Theorem \ref{th:0} is the possibility
to compare
the operators
$A$ and
$I(d^{2}/dx^{2})$.

\begin{corollary}[\cite{CM2006}]
There exists $k>0$ such that $A-kI(d^{2}/dx^{2})$ is
$L^{p}$-dissipative if and only if
$$
\underset{(x,{\lambda},{\omega})\in
(a,b)\times{\mathbb{C}}^{m}\times{\mathbb{C}}^{m}\atop
{|{\lambda}|=|{\omega}|=1}}
{\operatorname{ess\, inf}}
P(x,{\lambda},{\omega}) >0.
$$
There exists $k>0$ such that $kI(d^{2}/dx^{2})-A$ is
	 $L^{p}$-dissipative if and only if
	 $$
\underset{(x,{\lambda},{\omega})\in (a,b)\times{\mathbb{C}}^{m}\times{\mathbb{C}}^{m}
	 \atop
	 {|{\lambda}|=|{\omega}|=1}}
{\operatorname{ess\, sup}}
P(x,{\lambda},{\omega}) <\infty .
	 $$
There exists $k\in{\mathbb{R}}$ such that $A-kI(d^{2}/dx^{2})$ is
$L^{p}$-dissipative if and only if
$$
\underset{(x,{\lambda},{\omega})\in (a,b)\times{\mathbb{C}}^{m}\times{\mathbb{C}}^{m}
\atop
{|{\lambda}|=|{\omega}|=1}}
{\operatorname{ess\, inf}}
P(x,{\lambda},{\omega}) > -\infty.
$$
\end{corollary}

\subsection{Criteria in terms of eigenvalues
of ${\mathop{\mathscr A}\nolimits}(x)$}

If
the coefficients $a_{ij}$ of
the operator \eqref{defAord} are real,
it is possible to give
a necessary and sufficient condition for the $L^{p}$-dissipativity
of $A$ in terms of eigenvalues of the
matrix ${\mathop{\mathscr A}\nolimits}$.

\begin{theorem}[\cite{CM2006}]
Let ${\mathop{\mathscr A}\nolimits}$ be a
real matrix $\{a_{hk}\}$ with $h,k=1,\ldots,m$.
Suppose that ${\mathop{\mathscr A}\nolimits}={\mathop{\mathscr A}\nolimits}^{t}$ and ${\mathop{\mathscr A}\nolimits}\geq 0$
{\rm(}in the sense
that ${\langle}{\mathop{\mathscr A}\nolimits}(x)\xi,\xi\rangle\geq 0$
for almost every
$x\in(a,b)$ and for any $\xi\in{\mathbb{R}}^{m})$.
The operator $A$ is
$L^{p}$-dissipative if and only if
$$
\left(\frac{1}{2}-\frac{1}
	 {p}\right)^{2} ({\mu}_{1}(x)
	 +{\mu}_{m}(x))^{2} \leq {\mu}_{1}(x){\mu}_{m}(x)
$$
almost everywhere,
where ${\mu}_{1}(x)$ and ${\mu}_{m}(x)$ are the smallest and
largest
eigenvalues of the matrix ${\mathop{\mathscr A}\nolimits}(x)$ respectively.
 In the particular
case $m=2$, this condition is equivalent to
$$
\left(\frac{1}{2}-\frac{1}
	 {p}\right)^{2}({\operatorname{tr}} {\mathop{\mathscr A}\nolimits}(x))^{2} \leq \det {\mathop{\mathscr A}\nolimits}(x)
$$
almost everywhere.
\end{theorem}

\begin{corollary}[\cite{CM2006}]
Let ${\mathop{\mathscr A}\nolimits}$ be a real symmetric matrix.
Let ${\mu}_{1}(x)$ and ${\mu}_{m}(x)$ be the
	 smallest and
largest eigenvalues of ${\mathop{\mathscr A}\nolimits}(x)$ respectively.
	 There exists $k>0$ such that
	 $A-kI(d^{2}/dx^{2})$ is
	 $L^{p}$-dissipative if and only if
	\begin{equation}
\underset{x\in(a,b)}
{\operatorname{ess\, inf}}
		\left[(1+\sqrt{p\,p'}/2)\,{\mu}_{1}(x)+ (1-\sqrt{p\,p'}/2)\,
		{\mu}_{m}(x)
		\right]
			>0.
	 \label{primacondmu}
	\end{equation}
In the particular case $m=2$, the  condition \eqref{primacondmu} is
equivalent to
$$
\underset{x\in(a,b)}
{\operatorname{ess\, inf}}
\left[{\operatorname{tr}} {\mathop{\mathscr A}\nolimits}(x) -
	 \frac{\sqrt{p\,p'}}{2}\sqrt{({\operatorname{tr}}{\mathop{\mathscr A}\nolimits}(x))^{2}-
	 4\det{\mathop{\mathscr A}\nolimits}(x)} \right]
	 >0.
$$
\end{corollary}

Under an extra condition on the
matrix
${\mathop{\mathscr A}\nolimits}$, the following assertion holds.

	 \begin{corollary}[\cite{CM2006}]
Let ${\mathop{\mathscr A}\nolimits}$ be a real symmetric matrix.
Suppose that ${\mathop{\mathscr A}\nolimits}\geq 0$ almost everywhere.
Denote by ${\mu}_{1}(x)$ and ${\mu}_{m}(x)$ the
smallest and largest eigenvalues of
${\mathop{\mathscr A}\nolimits}(x)$ respectively.
If there exists $k>0$ such that
$A-kI(d^{2}/dx^{2})$ is
$L^{p}$-dissipative, then
\begin{equation}
\underset{x\in(a,b)}
{\operatorname{ess\, inf}}
\left[
{\mu}_{1}(x){\mu}_{m}(x) - \left(\frac{1}{2}-\frac{1}{p}\right)^{2}
({\mu}_{1}(x)+{\mu}_{m}(x))^{2}
\right]
>0.
\label{primacondmun}
\end{equation}
If, in addition, there exists
$C$ such that
\begin{equation}
{\langle}{\mathop{\mathscr A}\nolimits}(x)\xi,\xi{\rangle}\leq C|\xi|^{2}
\label{condC}
\end{equation}
for almost every $x\in(a,b)$ and for any $\xi\in{\mathbb{R}}^{m}$,
the converse assertion is also true.
In the particular case $m=2$, the  condition \eqref{primacondmun}
is equivalent to
$$
\underset{x\in(a,b)}
{\operatorname{ess\, inf}}
\left[\det {\mathop{\mathscr A}\nolimits}(x)
- \left(\frac{1}{2}-\frac{1}
{p}\right)^{2}
({\operatorname{tr}}{\mathop{\mathscr A}\nolimits}(x))^{2}\right] > 0 .
$$
\end{corollary}

We remark that $A-kI(d^{2}/dx^{2})$ is
$L^{p}$-dissipative means that $A$ is $p$-elliptic. 

Generally speaking, the assumption \eqref{condC} cannot be removed
even if ${\mathop{\mathscr A}\nolimits}\geq 0$.

\begin{example}
Consider $(a,b)=(1,\infty)$,
$m=2$,
${\mathop{\mathscr A}\nolimits}(x)=\{a_{ij}(x)\}$,
where
\begin{gather*}
a_{11}(x)=(1-2/\sqrt{pp'})x + x^{-1}, \quad
a_{12}(x)=a_{21}(x)=0,
\\[2pt]
a_{22}(x)=(1+2/\sqrt{pp'})x + x^{-1}.
\end{gather*}
Then
$${\mu}_{1}(x){\mu}_{2}(x) - \left(\frac{1}{2}-\frac{1}{p}\right)^{2}
({\mu}_{1}(x)+{\mu}_{2}(x))^{2} = (8+4x^{-2})/(p\,p')
$$
and \eqref{primacondmun} holds. But \eqref{primacondmu}
is not satisfied because
$$
(1+\sqrt{p\,p'}/2)\, {\mu}_{1}(x)+ (1-\sqrt{p\,p'}/2)\,
{\mu}_{2}(x) =2x^{-1} .
$$
\end{example}

\begin{corollary}[\cite{CM2006}]
Let ${\mathop{\mathscr A}\nolimits}$ be a real symmetric matrix.
Let
${\mu}_{1}(x)$ and ${\mu}_{m}(x)$ be the
smallest and largest eigenvalues of
${\mathop{\mathscr A}\nolimits}(x)$ respectively.
There exists $k>0$ such that $kI(d^{2}/dx^{2})-A$ is
$L^{p}$-dissipative if and only if
\begin{equation}
\underset{x\in(a,b)}
{\operatorname{ess\, sup}}
\left[(1-\sqrt{p\,p'}/2)\,{\mu}_{1}(x)+ (1+\sqrt{p\,p'}/2)\,
{\mu}_{m}(x) \right]
<\infty.
\label{secondacondmu}
\end{equation}
In the particular case $m=2$,
the  condition
\eqref{secondacondmu}
is equivalent to
$$
\underset{x\in(a,b)}
{\operatorname{ess\, sup}}
\left[{\operatorname{tr}} {\mathop{\mathscr A}\nolimits}(x) +
\frac{\sqrt{p\,p'}}{2}\sqrt{({\operatorname{tr}}{\mathop{\mathscr A}\nolimits}(x))^{2}-
4\det{\mathop{\mathscr A}\nolimits}(x)} \right]
<\infty.
		 $$
		 \end{corollary}

If ${\mathop{\mathscr A}\nolimits}$ is positive,
the following assertion holds.

\begin{corollary}[\cite{CM2006}]
Let ${\mathop{\mathscr A}\nolimits}$ be a real symmetric matrix.
Suppose that ${\mathop{\mathscr A}\nolimits}\geq 0$ almost everywhere.
Let
${\mu}_{1}(x)$ and ${\mu}_{m}(x)$ be the
smallest and largest eigenvalues
of ${\mathop{\mathscr A}\nolimits}(x)$ respectively.
There exists $k>0$ such that $kI(d^{2}/dx^{2})-A$ is
$L^{p}$-dissipative if and only if
$$
\underset{x\in(a,b)}
{\operatorname{ess\, sup}}
{\mu}_{m}(x)
<\infty.
	$$
\end{corollary}

\subsection{$L^{p}$-dissipativity of the operator \eqref{defspde}}

We represent necessary and
sufficient conditions for the $L^{p}$-dissipativity of
the system \eqref{defspde},
obtained in \cite{CM2006}.

Denote by $y_{h}$ the $(n-1)$-dimensional
vector $(x_{1},\ldots,x_{h-1}, x_{h+1}, \ldots ,x_{n})$
and set
${\omega}(y_{h})=\{ x_{h}\in{\mathbb{R}}\ |\ x\in{\Omega}\}$.

\begin{lemma}[\cite{CM2006}]
The operator \eqref{defspde} is $L^{p}$-dissipative if and only if
the ordinary differential operators
$$
A(y_{h})[u(x_{h})]=d({\mathop{\mathscr A}\nolimits}^{h}(x)du/dx_{h})/dx_{h}
$$
are $L^{p}$-dissipative in ${\omega}(y_{h})$
for almost every $y_{h}\in {\mathbb{R}}^{n-1}$
($h=1,\ldots,n$). This condition is void if
${\omega}(y_{h})=\emptyset$.
\end{lemma}

\begin{theorem}[\cite{CM2006}]
\label{thx:3}
The operator \eqref{defspde} is $L^{p}$-dissipative
if and only if
\begin{align}
\Re 
& {\langle}{\mathop{\mathscr A}\nolimits}^{h}(x_{0}) {\lambda},{\lambda}\rangle-(1-2/p)^{2}\Re 
\langle
{\mathop{\mathscr A}\nolimits}^{h}(x_{0}){\omega},{\omega}{\rangle}(\Re  \langle{\lambda},{\omega}\rangle)^{2}
\notag
\\[2pt]
&
 -(1-2/p)\Re ({\langle}{\mathop{\mathscr A}\nolimits}^{h}(x_{0}){\omega},{\lambda}{\rangle}-
{\langle}{\mathop{\mathscr A}\nolimits}^{h}(x_{0}){\lambda},{\omega}\rangle)
	 \Re  {\langle}{\lambda},{\omega}{\rangle} \geq 0
\label{79}
\end{align}
for almost every $x_{0}\in{\Omega}$ and for any ${\lambda},{\omega}\in{\mathbb{C}}^{m}$, $|{\omega}|=1$,
$h=1,\ldots,n$.
\end{theorem}

In the scalar case ($m=1$),
the operator \eqref{defspde} can be considered as
an operator from Section \ref{sec:scaldiv}.

In fact, if $Au=\sum_{h=1}^{n}{\partial}_{h}(a^{h}{\partial}_{h}u)$,
$a^{h}$ is a scalar function, then $A$ can be written in the
form \eqref{eq:nolower} with ${\mathop{\mathscr A}\nolimits}=\{c_{hk}\}$,
$c_{hh}=a^{h}$, $c_{hk}=0$ if $h\neq k$.
The conditions obtained in Section \ref{sec:scaldiv} can
be directly compared with
\eqref{79}.
We know that the operator $A$
is $L^{p}$-dissipative if and only if \eqref{eq:25} holds.
In this particular case, it is clear that
\eqref{eq:25} is equivalent to the following $n$ conditions:
\begin{equation}
\frac{4}{p\,p'}\, (\Re  a^{h})\, \xi^{2} + (\Re  a^{h})\, {\eta}^{2}
-2(1-2/p)(\Im  a^{h})\, \xi{\eta}\geq 0
\label{vecchiecond}
\end{equation}
almost everywhere and for any $\xi,{\eta}\in {\mathbb{R}}$, $h=1,\ldots,n$.
On the other hand, in this case, \eqref{79} reads as
\begin{align}
(\Re  a^{h})|{\lambda}|^{2}
&-(1-2/p)^{2}(\Re  a^{h})
(\Re  ({\lambda}\overline{{\omega}})^{2}
\notag\\[2pt]
&
-2(1-2/p)(\Im  a^{h})\Re ({\lambda}\overline{{\omega}})
\Im ({\lambda}\overline{{\omega}}) \geq 0
\label{nuovecond}
\end{align}
almost everywhere and for any
${\lambda},{\omega}\in{\mathbb{C}}$, $|{\omega}|=1$, $h=1,\ldots,n$. Setting
$\xi+i{\eta}={\lambda}\overline{{\omega}}$ and observing that
$|{\lambda}|^{2}=|{\lambda}\overline{{\omega}}|^{2}=(\Re ({\lambda}\overline{{\omega}}))^{2}+
(\Im ({\lambda}\overline{{\omega}}))^{2}$, we see that
the conditions \eqref{vecchiecond}
(hence
\eqref{eq:25}) are equivalent
to \eqref{nuovecond}.

If $A$ has real coefficients, we can characterize the
$L^{p}$-dissipativity in terms of the eigenvalues of
the matrices ${\mathop{\mathscr A}\nolimits}^{h}(x)$.

\begin{theorem}[\cite{CM2006}]
Let $A$ be an operator of the form
\eqref{defspde}, where
${\mathop{\mathscr A}\nolimits}^{h}$ are real
matrices $\{a^{h}_{ij}\}$ with $i,j=1,\ldots,m$.
Suppose that
${\mathop{\mathscr A}\nolimits}^{h}=({\mathop{\mathscr A}\nolimits}^{h})^{t}$ and ${\mathop{\mathscr A}\nolimits}^{h}\geq 0$
($h=1,\ldots,n$).
The operator $A$ is
$L^{p}$-dissipative if and only if
$$
\left(\frac{1}{2}-\frac{1}
	 {p}\right)^{2} ({\mu}_{1}^{h}(x)
	 +{\mu}_{m}^{h}(x))^{2} \leq
	 {\mu}_{1}^{h}(x)\,
	 {\mu}_{m}^{h}(x)
$$
for almost every $x\in{\Omega}$, $h=1,\ldots,n$,
where ${\mu}_{1}^{h}(x)$ and ${\mu}_{m}^{h}(x)$ are the
smallest and largest
eigenvalues of the matrix ${\mathop{\mathscr A}\nolimits}^{h}(x)$ respectively.
In the particular
case $m=2$, this condition is equivalent to
$$
\left(\frac{1}{2}-\frac{1}
{p}\right)^{2}({\operatorname{tr}} {\mathop{\mathscr A}\nolimits}^{h}(x))^{2} \leq \det
{\mathop{\mathscr A}\nolimits}^{h}(x)
$$
for almost every $x\in{\Omega}$, $h=1,\ldots,n$.
\end{theorem}

\section{The angle of dissipativity}\label{sec:angle}

By means of  the necessary and sufficient conditions we have obtained,
we can determine exactly the angle of dissipativity.
Determining the angle of dissipativity of an operator $A$,
means to find necessary and
sufficient conditions for the $L^{p}$-dissipativity of
the differential operator $zA$, where
$z\in{\mathbb{C}}$.

Consider first the scalar operator
$$
A={\nabla^{t}}({\mathop{\mathscr A}\nolimits}(x)\nabla),
$$
where ${\mathop{\mathscr A}\nolimits}(x)=\{a_{ij}(x)\}$ $(i,j=1,\ldots,n)$ is a matrix with complex
locally integrable entries defined in a domain
${\Omega}\subset{\mathbb{R}}^{n}$.
If ${\mathop{\mathscr A}\nolimits}$ is a
real matrix,
it is well known (cf.,
for example,
\cite{fattorini,fattorini2,okazawa})
that the dissipativity angle
is independent of
the operator and is given by
\begin{equation}
|\arg z| \leq \arctan \left(\frac{2\sqrt{p-1}}{|p-2|}\right).
\label{eq:fattoecc}
\end{equation}

If the entries of
the matrix ${\mathop{\mathscr A}\nolimits}$ are complex,
the situation
is different because the dissipativity angle
depends on the operator, as the next theorem shows.

\begin{theorem}[\cite{CM2006}]
\label{th:angdiss}
Let a
matrix ${\mathop{\mathscr A}\nolimits}$ be symmetric.
Suppose
that the operator $A$ is
$L^{p}$-dissipative. Let
$${\Lambda}_{1}=
\underset{(x,\xi)\in \Xi}
{\operatorname{ess\, inf}}
 \frac{{\langle}\Im  {\mathop{\mathscr A}\nolimits}(x)\xi,\xi\rangle}
{{\langle}\Re {\mathop{\mathscr A}\nolimits}(x)\xi,\xi\rangle}, \quad
{\Lambda}_{2}=
\underset{(x,\xi)\in\Xi}
{\operatorname{ess\, sup}}
 \frac{{\langle}\Im  {\mathop{\mathscr A}\nolimits}(x)\xi,\xi\rangle}
{{\langle}\Re {\mathop{\mathscr A}\nolimits}(x)\xi,\xi\rangle},
$$
where
$$
\Xi=\{ (x,\xi)\in {\Omega}\times{\mathbb{R}}^{n}\ |\
{\langle}\Re {\mathop{\mathscr A}\nolimits}(x)\xi,\xi\rangle
>0\}.
$$
The operator $zA$ is
$L^{p}$-dissipative if and only if
$$
{\vartheta}_{-}\leq \arg z \leq {\vartheta}_{+}\, ,
$$
where \footnote{Here, $0<{\mathop{\rm arccot}\nolimits} y <\pi$, ${\mathop{\rm arccot}\nolimits}(+\infty)=0$,
${\mathop{\rm arccot}\nolimits}(-\infty)=\pi$, and
$$
\underset{(x,\xi)\in \Xi}
{\operatorname{ess\, inf}}
 \frac{{\langle}\Im  {\mathop{\mathscr A}\nolimits}(x)\xi,\xi\rangle}
{{\langle}\Re {\mathop{\mathscr A}\nolimits}(x)\xi,\xi\rangle}=+\infty, \quad
\underset{(x,\xi)\in \Xi}
{\operatorname{ess\, sup}}
\frac{{\langle}\Im  {\mathop{\mathscr A}\nolimits}(x)\xi,\xi\rangle}
{{\langle}\Re {\mathop{\mathscr A}\nolimits}(x)\xi,\xi\rangle}=-\infty
$$
if $\Xi$ has zero measure.}
\begin{eqnarray*}
& \displaystyle {\vartheta}_{-}=
\begin{cases}
{\mathop{\rm arccot}\nolimits}\left( \frac{2\sqrt{p-1}}{|p-2|} -\frac{p^{2}}{|p-2|}
\, \frac{1}{2\sqrt{p-1}+|p-2|{\Lambda}_{1}}\right) -\pi
&\text{if}~~
p\neq 2,
\\[2pt]
{\mathop{\rm arccot}\nolimits}({\Lambda}_{1})-\pi
&\text{if}~~
p=2,
\end{cases}
&
\\[6pt]
& \displaystyle
{\vartheta}_{+}=\begin{cases}
{\mathop{\rm arccot}\nolimits} \left(-\frac{2\sqrt{p-1}}{|p-2|} + \frac{p^{2}}{|p-2|}
\,
\frac{1}{2\sqrt{p-1}-|p-2|{\Lambda}_{2}}\right)
&\text{if}~~
p\neq 2 \\[2pt]
{\mathop{\rm arccot}\nolimits}({\Lambda}_{2})
&\text{if}~~
p=2.
\end{cases}
&
\end{eqnarray*}
\end{theorem}

Note that
for a real matrix
${\mathop{\mathscr A}\nolimits}$
we have
${\Lambda}_{1}={\Lambda}_{2}=0$
and, consequently,
$$
\frac{2\sqrt{p-1}}{|p-2|} - \frac{p^{2}}{2\sqrt{p-1}|p-2|}=
- \frac{|p-2|}{2\sqrt{p-1}}\, .
$$
Theorem \ref{th:angdiss} asserts that $zA$
is dissipative if and only
if
$$
{\mathop{\rm arccot}\nolimits} \left(- \frac{|p-2|}{2\sqrt{p-1}}\right) -\pi \leq
\arg z \leq {\mathop{\rm arccot}\nolimits}\left(\frac{|p-2|}{2\sqrt{p-1}}\right),
$$
i.e., if and only if \eqref{eq:fattoecc} holds.

We can precisely determine the angle of dissipativity
also for the matrix ordinary
differential operator \eqref{defAord} with complex coefficients.

\begin{theorem}[\cite{CM2006}]
Let the operator \eqref{defAord} be $L^{p}$-dissipative.
The operator $zA$ is
	 $L^{p}$-dissipative if and only if
$$
	 {\vartheta}_{-}\leq \arg z \leq {\vartheta}_{+}
$$
	 where
\begin{align*}
{}&
{\vartheta}_{-}={\mathop{\rm arccot}\nolimits}\left(
\underset{(x,{\lambda},{\omega})\in\Xi}{\operatorname{ess\, inf}}
	 (Q(x,{\lambda},{\omega})/P(x,{\lambda},{\omega}))\right) - \pi,
\\[2pt]
&{\vartheta}_{+}={\mathop{\rm arccot}\nolimits}\left(
\underset{(x,{\lambda},{\omega})\in\Xi}{\operatorname{ess\, sup}}
	 (Q(x,{\lambda},{\omega})/P(x,{\lambda},{\omega}))\right),
\end{align*}
\begin{align*}
P(x,{\lambda},{\omega})
&=
	 \Re  {\langle}{\mathop{\mathscr A}\nolimits}(x) {\lambda},{\lambda}\rangle-(1-2/p)^{2}\Re {\langle}{\mathop{\mathscr A}\nolimits}(x) {\omega},{\omega}{\rangle}(\Re  \langle{\lambda},{\omega}\rangle)^{2}
		 \\[2pt]
		 &-
		 (1-2/p)\Re ({\langle}{\mathop{\mathscr A}\nolimits}(x) {\omega},{\lambda}{\rangle}-{\langle}{\mathop{\mathscr A}\nolimits}(x) {\lambda},{\omega}\rangle)
		 \Re  {\langle}{\lambda},{\omega}{\rangle},
	 \end{align*}
\begin{align*}
{}&Q(x,{\lambda},{\omega})=\Im  {\langle}{\mathop{\mathscr A}\nolimits}(x) {\lambda},{\lambda}\rangle-(1-2/p)^{2}\Im {\langle}{\mathop{\mathscr A}\nolimits}(x) {\omega},{\omega}\rangle
		 (\Re  \langle{\lambda},{\omega}\rangle)^{2}
\\[2pt]
&-
(1-2/p)\Im ({\langle}{\mathop{\mathscr A}\nolimits}(x) {\omega},{\lambda}{\rangle}-{\langle}{\mathop{\mathscr A}\nolimits}(x) {\lambda},{\omega}\rangle)
		 \Re  {\langle}{\lambda},{\omega}\rangle
\end{align*}
 and $\Xi$ is the set
	 $$
	 \Xi=\{ (x,{\lambda},{\omega})\in (a,b)\times {\mathbb{C}}^{m}\times {\mathbb{C}}^{m}\ |\
	 |{\omega}|=1,\ P^{2}(x,{\lambda},{\omega})+Q^{2}(x,{\lambda},{\omega})> 0\}.
	 $$
\end{theorem}

Finally Theorem \ref{thx:3} allows us
to determine the angle
of dissipativity of the operator \eqref{defspde}.

\begin{theorem}[\cite{CM2006}]
Let the operator  \eqref{defspde} be $L^{p}$-dissipative. The
operator $zA$ is $L^{p}$-dissipative if and only if
$$
{\vartheta}_{-}\leq \arg z \leq {\vartheta}_{+},
$$
 where
\begin{align*}
{}&
{\vartheta}_{-}=\max_{h=1,\ldots,n}\
{\mathop{\rm arccot}\nolimits}\left(
\underset{(x,{\lambda},{\omega})\in\Xi_{h}}
{\operatorname{ess\, inf}}
(Q_{h}(x,{\lambda},{\omega})/P_{h}(x,{\lambda},{\omega}))\right) - \pi,
\\[2pt]
&
{\vartheta}_{+}=\min_{h=1,\ldots,n}\
{\mathop{\rm arccot}\nolimits}\left(
\underset{(x,{\lambda},{\omega})\in\Xi_{h}}
{\operatorname{ess\, sup}}
(Q_{h}(x,{\lambda},{\omega})/P_{h}(x,{\lambda},{\omega}))\right)
\end{align*}
and
\begin{align*}
&P_{h}(x,{\lambda},{\omega})=
\Re  {\langle}{\mathop{\mathscr A}\nolimits}^{h}(x) {\lambda},{\lambda}\rangle-(1-2/p)^{2}\Re {\langle}{\mathop{\mathscr A}\nolimits}^{h}(x)
{\omega},{\omega}{\rangle}(\Re  \langle{\lambda},{\omega}\rangle)^{2}
\\[2pt]
&-
(1-2/p)\Re ({\langle}{\mathop{\mathscr A}\nolimits}^{h}(x){\omega},{\lambda}{\rangle}-{\langle}{\mathop{\mathscr A}\nolimits}^{h}(x){\lambda},{\omega}\rangle)
\Re  {\langle}{\lambda},{\omega}{\rangle},
\end{align*}
\begin{align*}
&
Q_{h}(x,{\lambda},{\omega})\!\!=\!\!\Im  {\langle}{\mathop{\mathscr A}\nolimits}^{h}(x) {\lambda},{\lambda}\rangle-(1-2/p)^{2}\!\!\!
\Im {\langle}{\mathop{\mathscr A}\nolimits}^{h}(x){\omega},{\omega}\rangle
(\Re  \langle{\lambda},{\omega}\rangle)^{2}
\\[2pt]
&
-
(1-2/p)\Im ({\langle}{\mathop{\mathscr A}\nolimits}^{h}(x){\omega},{\lambda}{\rangle}-{\langle}{\mathop{\mathscr A}\nolimits}^{h}(x){\lambda},{\omega}\rangle)
\Re  {\langle}{\lambda},{\omega}{\rangle},
\end{align*}
\begin{align*}
{}&
\Xi_{h}=
\{ (x,{\lambda},{\omega})\in \Omega
\times {\mathbb{C}}^{m}\times {\mathbb{C}}^{m}\ |\
|{\omega}|=1,\ P_{h}^{2}(x,{\lambda},{\omega})+Q_{h}^{2}(x,{\lambda},{\omega})> 0\}.
\end{align*}
\end{theorem}

\section{Maximum principles for linear elliptic equations
and systems}\label{sec:max}

As said in the Introduction, Kresin and Maz'ya have obtained
results on different forms of maximum principles for linear elliptic equations and systems. Here we recall some of their results. 

Let us consider the operator
\begin{equation}\label{eq:kresin1}
\mathfrak{A}_{0}\left(D_{x}\right)=\sum_{j, k=1}^{n} \mathcal{A}_{j k} \partial_{jk} 
\end{equation}
where $D_{x}=\left(\partial_{1}, \ldots, \partial_{n}\right)$ and $\mathcal{A}_{j k}=\mathcal{A}_{k j}$ are constant real $(m \times m)$-matrices. Assume that the operator $\mathfrak{A}_{0}$ is strongly elliptic, i.e., that for all $\zeta=\left(\zeta_{1}, \ldots, \zeta_{m}\right) \in \mathbb{R}^{m}$ and $\sigma=\left(\sigma_{1}, \ldots, \sigma_{n}\right) \in \mathbb{R}^{n}$, with $\zeta, \sigma \neq 0$, we have the inequality
$$
\left\lan \sum_{j, k=1}^{n} \mathcal{A}_{j k} \sigma_{j} \sigma_{k} {\zeta}, {\zeta}\right\ran >0
$$
Let $\Omega$ be a domain in $\mathbb{R}^{n}$ with boundary $\partial \Omega$ and closure $\overline{\Omega} .$ Let $\left[\mathrm{C}_{\mathrm{b}}(\overline{\Omega})\right]^{m}$ denote the space of bounded $m$-component vector-valued functions which are continuous in $\overline{\Omega}$. The norm on $\left[\mathrm{C}_{\mathrm{b}}(\overline{\Omega})\right]^{m}$ is $\|{u}\|=$ $\sup \{|{u}(x)|: x \in \overline{\Omega}\} .$ The notation $\left[\mathrm{C}_{\mathrm{b}}(\partial \Omega)\right]^{m}$ has a similar meaning. By $\left[\mathrm{C}^{2}(\Omega)\right]^{m}$ we denote the space of $m$-component vector-valued functions with continuous derivatives up to the second order in $\Omega$. 

Let
$$
\mathcal{K}(\Omega)=\sup \frac{\|{u}\|_{\left[\mathrm{C}_{\mathrm{b}}(\overline{\Omega})\right]^{m}}}{\|{u}\|_{\left[\mathrm{C}_{\mathrm{b}}(\partial \Omega)\right]^{m}}},
$$
where the supremum is taken over all vector-valued functions in the class $\left[\mathrm{C}_{\mathrm{b}}(\overline{\Omega})\right]^{m} \cap\left[\mathrm{C}^{2}(\Omega)\right]^{m}$ satisfying the system $\mathfrak{A}_{0}\left(D_{x}\right) {u}={0}$.

Clearly, $\mathcal{K}(\Omega)$ is the best constant in the inequality
$$
|{u}(x)| \leq \mathcal{K}(\Omega) \sup \{|{u}(y)|: y \in \partial \Omega\}
$$
where $x \in \Omega$ and ${u}$ is a solution of the system $\mathfrak{A}_{0}\left(D_{x}\right) {u}=0$ in the class $\left[\mathrm{C}_{\mathrm{b}}(\overline{\Omega})\right]^{m} \cap\left[\mathrm{C}^{2}(\Omega)\right]^{m}$

If $\mathcal{K}(\Omega)=1$, then the classical maximum modulus principle holds for the system $\mathfrak{A}_{0}\left(D_{x}\right) {u}=\mathbf{0}$.

Kresin and Maz'ya proved the following 
criterion for the validity of this classical modulus principle.
\begin{theorem}
    Let $\Omega$ be a domain in $\mathbb{R}^{n}$ with compact closure and $\mathrm{C}^{1}$-boundary. The equality $\mathcal{K}(\Omega)=1$ holds if and only if the operator $\mathfrak{A}_{0}\left(D_{x}\right)$ is defined by 
\begin{equation}\label{eq:(2.10)}
\mathfrak{A}_{0}\left(D_{x}\right)=\mathcal{A} \sum_{j, k=1}^{n} a_{j k} \partial_{jk}
\end{equation}
where $\mathcal{A}$ and $\{a_{j k}\}$ are positive-definite constant matrices of orders $m$ and $n$, respectively.
\end{theorem}

Suppose now that the operator \eqref{eq:kresin1} has complex
coefficients, i.e., suppose that $\mathcal{A}_{j k}=\mathcal{A}_{k j}$ are constant  complex $(m \times m)$-matrices. Assume that the operator
is strongly elliptic. This means that
$$
\Re \left\lan \sum_{j, k=1}^{n} \mathcal{A}_{j k} \sigma_{j} \sigma_{k} {\zeta}, {\zeta}\right\ran >0
$$
for all $\zeta=\left(\zeta_{1}, \ldots, \zeta_{m}\right) \in \mathbb{C}^{m}$ and $\sigma=\left(\sigma_{1}, \ldots, \sigma_{n}\right) \in \mathbb{R}^{n}$, with $\zeta, \sigma \neq 0$

A necessary and sufficient condition for validity of the classical modulus principle for operator \eqref{eq:kresin1} with complex coefficients in a bounded domain runs as follows.

\begin{theorem}
    Let $\Omega$ be a domain in $\mathbb{R}^{n}$ with compact closure and $\mathrm{C}^{1}$-boundary. The equality $\mathcal{K}(\Omega)=1$ holds if and only if the operator $\mathfrak{A}_{0}\left(D_{x}\right)$ has the form \eqref{eq:(2.10)}, where now 
    $\mathcal{A}$ is a constant complex-valued $({m} \times{m})$-matrix such that $\Re(\mathcal{A} \zeta, \boldsymbol{\zeta})>$ 0 for all $\zeta \in \mathbb{C}^{m}, \zeta \neq \mathbf{0}$, and $\{a_{j k}\}$ is a real positive-definite $(n \times n)$ matrix.
\end{theorem}

These results have been extended to more general systems and we refer to the survey \cite{kresmazsurv} for all the details.

\section{Other results}\label{sec:others}

In this section we briefly mention other results we have obtained.

In \cite{CM2018} we found necessary and sufficient conditions
for the $L^{p}$-dissipativity of systems of the first order.
 Namely we have considered the matrix operator
 \begin{equation}\label{eq:opB}
Eu=\B^{h}(x)\de_{h}u + \Dm(x) u\, ,
\end{equation}
where $\B^h(x)=\{ b^h_{ij}(x)\}$  and 
$\Dm(x)=\{ d_{ij}(x)\}$ are matrices with complex locally
integrable entries defined in the domain  $\Om$
of $\R^n$ and $u=(u_{1},\ldots,u_{m})$ ($1\leq i,j \leq m,\ 1\leq h 
\leq n$). It states that, if $p\neq 2$, $E$ is $L^{p}$-dissipative if, and only if, 
\begin{equation}
    \B^{h}(x)=b_{h}(x) I \  \text{a.e.},
    \label{eq:0}
\end{equation} 
$b_{h}(x)$ being real locally integrable 
functions, and
the inequality
$$
\Re \lan (p^{-1} \de_{h}\B^{h}(x) - \Dm(x)) \ze,\ze\ran  \geq 0
$$
holds for any $\ze\in \CC^{m}$, $|\ze|=1$ and for almost any $x\in\Om$.
If $p=2$ condition \eqref{eq:0} is replaced
by the more general requirement that the matrices $\B^{h}(x)$ are  
self-adjoint a.e..

We have applied this result also to second order operators, obtaining a sufficient
condition for their $L^p$-dissipativity. We have also determined the angle
of dissipativity of operator \eqref{eq:opB}.

In \cite{CM2019} we have considered the ``complex oblique derivative''
operator
\begin{equation}\label{eq:cod}
\lambda\cdot\nabla u=
\frac{\partial u}{\partial x_{n}} +
\sum_{j=1}^{n-1}a_{j}\frac{\partial u}{\partial x_{j}}
\end{equation}
where $\lambda=(1,a_1,\ldots,a_{n-1})$ and $a_j$ are  complex valued 
functions. We gave necessary and, separately, sufficient 
conditions under which such boundary operator is $L^p$-dissipative
on $\R^{n-1}$.
If the coefficients $a_j$ are real valued, 
we have obtained a  necessary and sufficient condition:
 the operator \eqref{eq:cod}  is $L^{p}$-dissipative if and only if
there exists a real vector $\Gamma\in L^2_{\text{loc}}(\R^n)$ such that
$$
    -\partial_{j}(\Re a_{j}) \, \delta(x_n) \leq \frac{2}{p'}( \dive \Gamma - 
    |\Gamma|^{2})
$$
in the sense of distributions.

In the same paper we have considered also the class of integral
operators which can be
written as
\begin{equation}\label{eq:intoper}
\int_{\R^{n}}^*[u(x)-u(y)]\, K(dx,dy)
\end{equation}
where the integral has to be understood as a principal value in the sense of Cauchy
and the kernel $K(dx,dy)$ is a Borel positive measure defined 
on $\R^n\times\R^n$ satisfying certain conditions. The class
of operators we considered includes the fractional powers of Laplacian
$(-\Delta)^s$, with $0<s<1$. 
For the latter we previously had proved  the following theorem 
 \begin{theorem}
[\cite{CMbook}, p.230--231]
 \label{th:tipokato}
   Let $0<\al<1$. We have,  for any $u\in \Cspt^\infty(\R^n)$,
   $$
   \int_{\R^n} \lan (-\D)^\al u, u\ran |u|^{p-2} dx
   \geq  \frac{2\, c_\al}{p\, p'} \, \Vert |u|^{p/2}\Vert^2_{{\mathcal L}^{\al,2}(\R^n)}\, ,
   $$
  where
 $$
    c_\al = -\pi^{-n/2} 4^{\al} \GA(\al+n/2)/ \GA(-\al) >0
$$
and
  $$
 \Vert  v \Vert_{{\mathcal L}^{\al,2}(\R^n)}
  =  \left(\iint_{\R^n\times \R^n} |v(y)-v(x)|^2 \frac{dxdy}{|y-x|^{n+2\al}}\right)^{1/2}.
 $$
    \end{theorem}

In \cite{CM2019} we have established the $L^p$-positivity
of operator \eqref{eq:intoper}, extending in this way  Theorem \ref{th:tipokato}.

\section{The functional dissipativity}\label{sec:funct}

In \cite{CM2021} we have introduced the new concept of functional
dissipativity. Roughly speaking the idea is to replace
 $|u|^{p-2}$ by a more general  $\vf(|u|)$, $\vf$ being a positive function. 

Let us consider the operator \eqref{eq:nolower} with $L^\infty$ complex
valued coefficients. 
We say that it
is functional dissipative or $L^\Phi$-dissipative if 
$$
\Re \int_\Om \lan \A \nabla u, \nabla(\vf(|u|)\, u)\ran\, dx \geq 0
$$
for any $u\in \Hspt^{1}(\Om)$ such that $\vf(|u|)\, u \in \Hspt^{1}(\Om)$.
Here $\vf$ is a positive function defined on $\R^{+}=(0,+\infty)$ which  satisfies the following conditions:

\renewcommand{\labelenumi}{(\roman{enumi})}
\renewcommand{\theenumi}{(\roman{enumi})}
\begin{enumerate}\label{items}
	\item\label{item1} $\vf \in C^{1}((0,+\infty))$;
	\item $(s\, \vf(s))'>0$ for any $s>0$;
	\item the range of the strictly increasing function $s\, \vf(s)$ is  $(0,+\infty)$;	
	\item\label{item4} there exist two positive constants $C_{1}, C_{2}$  and a real number $r>-1$ such that
$$
C_{1} s^{r}\leq (s\vf(s))' \leq C_{2}\, s^{r}, \qquad s\in (0,s_{0})
$$
for a certain $s_{0}>0$. If $r=0$ we require more restrictive 
conditions: there exists the finite limit $\lim_{s\to 
0^+}\vf(s)=\vf_{+}(0)>0$ 
and  $\lim_{s\to 0^+}s\, \vf'(s)=0$.
\item\label{item5} 
There exists $s_{1}>s_{0}$ such that 
$$
   \vf'(s)\geq 0  \text{ or }  \vf'(s)\leq 0 \qquad \forall\ s\geq s_{1}  .
$$
\end{enumerate}

The reason for requiring that function $s\, \vf(s)$ is increasing is that
in such a way  the function 
\begin{equation}\label{eq:Phi}
\Phi(s)=\int_{0}^{s} \si\, \vf(\si)\, d\si
\end{equation}
is a Young function (i.e., a convex positive function such that $\Phi(0)=0$ and 
$\Phi(+\infty)=+\infty$) . We note that, if $t\, \psi(t)$ is the inverse function of $s\, \vf(s)$, then
$$
\Psi(s)= \int_{0}^{s} \si\, \psi(\si)\, d\si
$$
is the conjugate Young function of $\Phi$.

The condition \ref{item4} prescribes the behaviour of the function $\vf$ in a neighborhood of the origin,
while \ref{item5} concerns the behaviour for large $s$. 

The function $\vf(s)=s^{p-2}$ ($p>1$) provides an example of such a function.

A motivation for the study of the concept 
of functional dissipativity comes from the decrease of the  Luxemburg 
norm of solutions of the Cauchy–Dirichlet problem
\begin{equation}\label{eq:cauchy}
 \begin{cases}
u'=Au \\ 
u(0)=u_{0} \, . 
\end{cases}
\end{equation}

 Indeed let us
consider  the Orlicz  space of functions $u$ for which
there exists $\al>0$ such that
$$
\int_{\Om}\Phi(\al\, |u|)\, dx < +\infty\, .
$$

For the general theory of Orlicz spaces we refer to Krasnosel'ski\u{\i}, Ruticki\u{\i} \cite{krasno}
and Rao, Ren \cite{rao}.
As in \eqref{eq:dt}, if $u(x,y)$ is a solution of
the Cauchy-Dirichlet problem \eqref{eq:cauchy}, we have the decrease of the 
integrals 
$$
\int_{\Om}\Phi( |u(x,t)|)\, dx 
$$
 if
$$
\Re \int_{\Om} \lan Eu,u\ran |u|^{-1}\Phi'( |u|)\, dx \leq 0.
$$

This implies  the decrease of the
Luxemburg norm in the related Orlicz space
$$
\Vert u(\cdot, 
    t)\Vert =
\inf \left\{ \lan>0\ |\ \int_{\Om}\Phi(|u(x,t)|/\la)\, dx \leq 1 
\right\}.
$$

In paper \cite{CM2021} we proved the following technical lemma, which
played a key role.

\begin{lemma}\label{lemma:fond}
   The operator $A$  is $L^{\Phi}$-dissipative if and only if
\begin{gather*}
\Re \int_{\Om}\Big[ \lan\A\nabla v, \nabla v \ran +
\Lambda(|v|)\, \lan(\A-\A^{*})\nabla |v|, |v|^{-1}\overline{v} \nabla v ) 
\ran +\\
-\Lambda^{2}(|v|)\, \lan\A\nabla |v|, \nabla |v| \ran \Big] dx \geq 0,
\qquad \forall v \in \Hspt^{1}(\Om),
\end{gather*}
where the function $\Lambda$ is is the function defined
by the relation
$$
\Lambda\left(s\sqrt{\vf(s)}\right)= - \frac{s\, \vf'(s)}{s\,\vf'(s)+2\, 
\vf(s)}\, .
$$
\end{lemma}

We remark that if $\vf(s)=s^{p-2}$, the function $\Lambda$ is constant and  $$\Lambda(t)=-(1-2/p), \quad
1-\Lambda^{2}(t)=4/(p\, p').
$$

As Corollaries of  Lemma \ref{lemma:fond} we have obtained necessary and, 
separately,  sufficient
conditions for the functional dissipativity of the operator $E$.

\begin{corollary}
    If the operator $A$  is $L^{\Phi}$-dissipative, we have
	\begin{equation}\label{eq:E>=0}
\lan \Re \A (x) \xi,\xi\ran \geq 0
\end{equation}
for almost every $x\in\Om$ and for any  $\xi\in\R^{n}$.
\end{corollary}

\begin{corollary}
If
\begin{equation}\label{eq:polsuf}
\begin{gathered}
{}[1-\Lambda^{2}(t)] \lan\Re \A(x)\, \xi,\xi\ran 
+ \lan\Re \A(x)\, \eta,\eta\ran +\\
[1+\Lambda(t)] \lan \Im \A(x)\, \xi,\eta\ran 
+
[1-\Lambda(t)] \lan \Im \A^{*}(x)\, \xi,\eta\ran
\geq 0
\end{gathered}
\end{equation}
for almost every $x\in\Om$ and for any $t>0, \xi, \eta \in \R^{n}$, 
the operator $A$  is $L^{\Phi}$-dissipative.
\end{corollary}

\begin{corollary}
If the operator $A$ has real coefficients
and satisfies condition \eqref{eq:E>=0} for almost every $x\in\Om$ and for any  $\xi\in\R^{n}$,
than it is $L^{\Phi}$-dissipative with respect to any $\Phi$.
\end{corollary}

The main result obtained in \cite{CM2021} is the following necessary and sufficient condition

\begin{theorem}\label{th:mainfunct}
  Let the matrix $\Im \A$ be symmetric, i.e.,
  $\Im \A^{t}=\Im \A$. 
  Then the operator $A$  is 
  $L^{\Phi}$-dissipative if, and only if, 
  \begin{equation}\label{eq:cmcond}
|s\, \vf'(s)| \, | \lan\Im \A (x)\, \xi,\xi\ran | 
\leq 2\, \sqrt{\vf(s)\, [s\, \vf(s)]'} 
\, \lan \Re \A(x) \, \xi,\xi\ran 
\end{equation}
for almost every $x\in\Om$ and for any  $s>0, \xi\in\R^{n}$.
\end{theorem}

Suppose that the condition $\Im\A=\Im\A^t$ is not satisfied.
Arguing as in the proof of Theorem \ref{th:mainfunct}, one can prove that
condition \eqref{eq:cmcond} is still necessary for the $L^{\Phi}$-dissipativity of the operator $E$. 
However in general it is not sufficient, whatever the function $\vf$ may be.
This is shown by the next example. 

\begin{example}
Let $n=2$ , $\Om$ be a bounded domain, $\la$ be a real parameter and
$$
	\A=\left(\begin{array}{cc}
	1 & i\la x_1 \\ -i\la x_1 & 1
	\end{array}\right)
	$$

	Since $\lan \Re\A \xi,\xi\ran = |\xi|^{2}$ and 
	$\lan \Im\A \xi,\xi\ran = 0$ for any $\xi\in\R^{n}$, condition
	\eqref{eq:cmcond} is  satisfied.
	
If the corresponding operator $Eu=\Delta u + i\, \la \, \partial_2 u$ 
is $L^{\Phi}$-dissipative,
then
\begin{equation}\label{eq:exfinale}
\Re \int_\Om \lan \Delta u + i\, \la \, \partial_2 u, u\ran\, \vf(|u|)\, 
dx \leq 0, \qquad \forall\ u\in \Cspt^\infty(\Om).
\end{equation}

Take $u(x)=\ro(x)\, e^{i\, t \, x_{2}}$, where $\ro\in 
\Cspt^\infty(\Om)$ is real valued and $t\in\R$.
Since $ \lan Eu,u\ran = \ro[\Delta \ro + 2\, i\, t\, \partial_{2}\ro - 
t^{2}\ro
+ i\,\la\,(\partial_{2}\ro + i t \ro)]$,
condition \eqref{eq:exfinale} implies
\begin{equation}\label{eq:imposs}
\int_{\Om}\ro\, \Delta \ro\, \vf(|\ro|)\, dx - \la\, t \int_{\Om}\ro^{2} 
\vf(|\ro|)\, dx -t^{2}\int_{\Om}\ro^{2} 
\vf(|\ro|)\, dx \leq 0
\end{equation}
for any $t, \la\in\R$. The function $\vf$ being positive, we can 
choose $\ro$ in such a way
$$
\int_{\Om}\ro^{2} 
\vf(|\ro|)\, dx >0.
$$

Taking  
$$
\la^{2}>4 \int_{\Om}\ro\, \Delta \ro\, \vf(|\ro|)\, dx 
\left(\int_{\Om}\ro^{2} 
\vf(|\ro|)\, dx\right)^{-1},
$$
inequality \eqref{eq:imposs} is impossible for all $t\in\R$. Thus
$E$ is not $L^{\Phi}$-dissipative, although \eqref{eq:cmcond} is satisfied.
\end{example}

We have also
\begin{corollary}\label{co:4}
	Let the matrix $\Im \A$ be symmetric, i.e.,
  $\Im \A^{t}=\Im \A$. 
      If 
\begin{equation}\label{eq:lam0}
	  \la_{0}=
	\sup_{s>0} 
	\frac{|s\, \vf'(s)|}
 	{2\, \sqrt{\vf(s)\, [s\, \vf(s)]'}} < +\infty, 
\end{equation}
then the operator $E$   is 
  $L^{\Phi}$-dissipative if, and only if, 
    \begin{equation}\label{eq:cmcond2}
\la_{0}\,  | \lan\Im \A (x)\, \xi,\xi\ran | 
\leq 
 \lan \Re \A(x) \, \xi,\xi\ran 
\end{equation}
for almost every $x\in\Om$ and for any $ \xi\in\R^{n}$. If $\la_{0}=+\infty$  the operator $E$  is $L^{\Phi}$-dissipative
 if and only if $\Im\A= 0$ and condition \eqref{eq:E>=0}
is satisfied.
\end{corollary}

If we use the function $\Phi$ (see \eqref{eq:Phi}), 
condition \eqref{eq:cmcond} can be written as
$$
|s\, \Phi''(s) - \Phi'(s)| \, | \lan\Im \A (x)\, \xi,\xi\ran | 
\leq 2 \sqrt{s\, \Phi'(s)\, \Phi''(s)} 
\, \lan \Re \A(x) \, \xi,\xi\ran 
$$
for almost every $x\in\Om$ and for any $s>0,  \xi\in\R^{n}$. In the same way, 
formula  \eqref{eq:lam0} becomes
$$	  \la_{0}=
	\sup_{s>0} 
	\frac{|s\, \Phi''(s) - \Phi'(s)|}
 	{2 \sqrt{s\, \Phi'(s)\, \Phi''(s)}} < +\infty. $$

We consider now some examples in which we indicate both 
the functions $\Phi$ and $\vf$. It is easy to verify that in each example the
function $\vf$ satisfies conditions  \ref{item1}-\ref{item5}
(see p.\pageref{items}).

\begin{example}
	If $\Phi(s)=s^p$, i.e., $\vf(s)=p\, s^{p-2}$, which corresponds to 
	$L^{p}$ norm, the function in \eqref{eq:lam0} is constant and
	$\la_ {0}=|p-2|/(2\sqrt{p-1})$. In this way we reobtain  Theorem \ref{th:main}.
\end{example}

\begin{example}
	Let us consider $\Phi(s)=s^{p}\log(s+e)$ ($p>1$), which is the
	Young function corresponding to the Zygmund space $L^{p}\, \log L$.
	This is equivalent to say $\vf(s)=p s^{p-2}\log(s+e) + 
	s^{p-1}(s+e)^{-1}$.
	By a direct computation we find
	\begin{equation}
\begin{gathered}
\frac{|s\, \Phi''(s) - \Phi'(s)|}
 	{2\sqrt{s\, \Phi'(s)\, \Phi''(s)}} =\\
\frac{\left| p(p-2)\log(s+e)+\frac{(2p-1)s}{s+e} - \frac{s^2}{(s+e)^{2}} \right| }
{2  \sqrt{\left(p\log(s+e) + \frac{s}{s+e}\right) \left(p(p-1)\log(s+e) + \frac{2ps}{s+e} - \frac{s^2}{(s+e)^{2}}\right) }} \, .
\label{eq:supPhi}
\end{gathered}
\end{equation}

Since
$$
\lim_{s\to 0^+}\frac{|s\, \Phi''(s) - \Phi'(s)|}
 	{2\sqrt{s\, \Phi'(s)\, \Phi''(s)}} = \lim_{s\to +\infty}\frac{|s\, \Phi''(s) - \Phi'(s)|}
 	{2\sqrt{s\, \Phi'(s)\, \Phi''(s)}} = \frac{|p-2|}{2\sqrt{p-1}}
$$
the function is bounded. Then we have  the 
$L^{\Phi}$-dissipativity of the operator $A$  
if, and only if, \eqref{eq:cmcond2} holds,
where  $\la_ {0}$ is the $\sup$ of the function  \eqref{eq:supPhi}
in $\R^{+}$.\end{example}

%

\begin{example}
Let us consider the function $\Phi(s)=\exp(s^{p})-1$,
i.e., $\vf(s)=p\, s^{p-2} \exp(s^{p})$.
In this case
	$$\frac{|s\, \Phi''(s) - \Phi'(s)|}
 	{2\sqrt{s\, \Phi'(s)\, \Phi''(s)}} =
	\frac{|p\, s^p+ p -2|}{2 \sqrt{(p\, s^p +p -1)}}$$
and $\la_ {0}=+\infty$. In view of Corollary \eqref{co:4}, the operator $A$  is
$L^{\Phi}$-dissipative, i.e.,
$$
	\Re \int_{\Om} \lan \A\nabla u, \nabla[u\, |u|^{p-2} \exp(|u|^p) ]\ran dx \geq 0
	$$
 for any $u\in\Hspt^{1}(\Om)$ such that $|u|^{p-2} \exp(|u|^p)\, u\in\Hspt^{1}(\Om)$, if and only if the operator $A$ has real 
coefficients and condition \eqref{eq:E>=0}
is satisfied.
\end{example}

\begin{example}
	Let $\Phi(s)=s-\arctan s$, i.e., $\vf(s)=s/(s^{2}+1)$. In this case
	$$
	\frac{|s\, \Phi''(s) - \Phi'(s)|}
 	{2\sqrt{s\, \Phi'(s)\, \Phi''(s)}} =
	\frac{|s^{2}-1|}{2 \sqrt{2(s^{2}+1})}
	$$
	and $\la_ {0}=+\infty$. As in the previous example, we have that
	$$
	\Re \int_{\Om} \lan\A\nabla u , \nabla\left(\frac{|u|\, u}{|u|^{2}+1}\right)\ran dx \geq 0
	$$
for any $u\in\Hspt^{1}(\Om)$ such that $|u|\, u /(|u|^{2}+1)\in\Hspt^{1}(\Om)$,  if and only if the operator $A$ has real 
coefficients and condition \eqref{eq:E>=0}
is satisfied.
\end{example}

\begin{example}
	Let $\Phi(s)=s^{4}/(s^{2}+1)$, i.e., $\vf(s)=2\,s^{2}(2+s^{2})/(s^{2}+1)^{2}$. In this case
	$$
	\frac{|s\, \Phi''(s) - \Phi'(s)|}
 	{2\sqrt{s\, \Phi'(s)\, \Phi''(s)}} =
	\frac{2}{\sqrt{(s^{2}+1)(s^{2}+2)(s^{4}+3s^{2}+6)}}\, .
	$$
	This function is decreasing and  $\la_ {0}$ is equal to its 
	value at $0$, i.e., $\la_ {0}=1/\sqrt{3}$. 
	The operator $A$  is
$L^{\Phi}$-dissipative, i.e.,
	$$
	\Re \int_{\Om} \lan \A\nabla u, \nabla\left(\, \frac{|u|^{2}(2+|u|^2)u}{(|u|^{2}+1)^{2}}\right)\ran dx \geq 0
	$$ 
	for any $u\in\Hspt^{1}(\Om)$ such that $|u|^{2}(2+|u|^2)u /(|u|^{2}+1)^{2} \in\Hspt^{1}(\Om)$, if and only 
	if
	$$
	| \lan\Im \A (x)\, \xi,\xi\ran | 
\leq \sqrt{3} 
\, \lan \Re \A(x) \, \xi,\xi\ran 
$$
for almost any $x\in\Om$ and for any $\xi\in\R^{n}$.
\end{example}

\begin{example}
	Let $\Phi(s)=s^{2}(s^{2}+2)/(s^{2}+1)-2 \log(s^{2}+1)$, 
	i.e., $\vf(s)=2\,s^{4}/(s^{2}+1)^{2}$. In this case
	$$
	\frac{|s\, \Phi''(s) - \Phi'(s)|}
 	{2\sqrt{s\, \Phi'(s)\, \Phi''(s)}} =
	\frac{2}{\sqrt{(s^{2}+1)(s^{2}+5)}}\, .
	$$
	This function is decreasing and  $\la_ {0}$ is equal to its 
	value at $0$, i.e., $\la_ {0}=2/\sqrt{5}$. 
	The operator $A$  is
$L^{\Phi}$-dissipative, i.e.,
	$$
	\Re \int_{\Om} \lan \A\nabla u, \nabla\left(\frac{|u|^{4}u}{(|u|^{2}+1)^{2}}\right)\ran dx \geq 0
	$$  
	for any $u\in\Hspt^{1}(\Om)$ such that $|u|^{4}u /(|u|^{2}+1)^{2} \in\Hspt^{1}(\Om)$,
	if and only 
	if
	$$
	2\, | \lan\Im \A (x)\, \xi,\xi\ran | 
\leq \sqrt{5} 
\, \lan \Re \A(x) \, \xi,\xi\ran 
$$
for almost any $x\in\Om$ and for any $\xi\in\R^{n}$.
\end{example}

By analogy to the $L^p$ case, if we have an operator
with lower order term \eqref{eq:lowerterms} and if the principal part  is such that
the left-hand side of \eqref{eq:polsuf} is not merely non negative but strictly positive, i.e. 
\begin{gather*}
{}[1-\Lambda^{2}(t)] \lan\Re \A(x)\, \xi,\xi\ran 
+ \lan\Re \A(x)\, \eta,\eta\ran +\\
[1+\Lambda(t)] \lan \Im \A(x)\, \xi,\eta\ran 
+
[1-\Lambda(t)] \lan \Im \A^{*}(x)\, \xi,\eta\ran
\geq \kappa (|\xi|^2 + |\eta|^2)
\end{gather*}
for a certain $\kappa>0$ and for almost every $x\in\Om$ and for any $t>0,  \xi, \eta \in \R^{n}$,
we say that the operator $A$ is (strongly) $\Phi$-elliptic.

We note that, if  $A$ is a   (strongly) $\Phi$-elliptic operator, then there exists a constant $\kappa$ such that
    for any nonnegative $\chi\in L^{\infty}(\Om)$ and any complex valued $u\in H^1(\Om)$ such that
   $\vf(|u|)\, u \in H^1(\Om)$ we have
      $$
\Re \int_{\Om} \lan \A \nabla u, \nabla(\vf(|u|)\, u)\ran\, \chi(x) dx \geq \kappa \int_{\Om}
| \nabla(\sqrt{\vf(|u|)}\, u)|^{2} \chi(x)\, dx
$$
(see \cite[Corollary 4]{CM2021}).

\section{Concluding remarks}\label{sec:conclude}

Our condition \eqref{intro:form} and its strengthened variant
are getting more and more important in many respects.
We said  already something about $p$-ellipticity,
but there are also other applications.

We mention that 
{H\"omberg, Krumbiegel and Rehberg \cite{HKR} used some of the techniques introduced in \cite{CM2005} to show 
the $L^{p}$-dissipativity of a certain operator connected to the problem
of the existence of an optimal control for the heat equation with dynamic boundary condition.


Beyn and  Otten  \cite{BO2016,BO2018} considered the semilinear system
$$
A \Delta v(x)  +\lan Sx, \nabla v(x)\ran + f(v(x))=0, \qquad  x\in R^N,
$$
where $A$ is a $m\times m$ matrix, $S$ is a $N\times N$ skew-symmetric matrix
and $f$ is a sufficiently smooth vector function. Among the assumptions they made, they require
the existence of a constant $\ga_A>0$ such that
$$
|z|^2 \Re \lan w, Aw\ran + (p-2) \Re \lan w,z\ran \Re \lan z, Aw\ran \geq \ga_A |z|^{2} |w|^{2}
$$
for any $z,w\in \C^m$. This condition originates from our \eqref{79}.

The results of \cite{CM2005}
 allowed  Nittka \cite{nittka} to consider 
the case of partial differential operators with complex coefficients. 

Ostermann and  Schratz \cite{oster} have obtained the stability  of a numerical procedure for
solving a certain evolution problem.  
The necessary and sufficient condition \eqref{eq:24}  show that
their  result does not require  the contractivity of the corresponding semigroup.
  
 Chill,  Meinlschmidt and Rehberg \cite{chill} used some ideas from \cite{CM2005}
 in the study of the numerical range of second order elliptic operators with mixed boundary conditions in $L^p$.

ter Elst, Haller-Dintelmann, Rehberg and Tolksdorf \cite{EHR}
considered second
 order divergence form operators with complex coefficients, complemented with Dirichlet, Neumann or mixed boundary conditions.  They proved several results related to the generation of  strongly continuous semigroup on $L^p$.

 \section*{Funding}

The second author has been supported by the RUDN University Strategic Academic
Leadership Program.

\end{document}